\documentclass[reqno,a4paper]{amsart}
\usepackage{amssymb}%
\usepackage[hmarginratio=1:1]{geometry}%
\usepackage{ifpdf}
\usepackage{listings}
\usepackage{color,xcolor}
\usepackage{colortbl}
\usepackage{multirow}
\usepackage{booktabs}
\usepackage{float}
\usepackage{graphicx}
\usepackage{algorithmic}
\usepackage[linesnumbered,ruled,vlined]{algorithm2e}
\usepackage{subfig}
\usepackage{epstopdf}
\usepackage[colorlinks,linkcolor=red,anchorcolor=red,citecolor=red]{hyperref}
\usepackage[numbers,sort&compress]{natbib}
\usepackage{enumitem}

\theoremstyle{plain}
\newtheorem{thm}{Theorem}[section]

\newtheorem{cor}{Corollary}[section]

\theoremstyle{remark}
\newtheorem{remark}{Remark}[section]
\theoremstyle{definition}

\numberwithin{equation}{section}
\allowdisplaybreaks[4]

\begin{document}
\title[]{Recurrence Relations for the Maclaurin Coefficients of Products of Elementary Functions and Hypergeometric Functions}
\author{Zhong-Xuan Mao, Jing-Feng Tian*}

\address{Zhong-Xuan Mao \\
Department of Mathematics and Physics\\
North China Electric Power University \\
Yonghua Street 619, 071003, Baoding, P. R. China}
\email{maozhongxuan000\symbol{64}gmail.com}

\address{Jing-Feng Tian\\
Department of Mathematics and Physics\\
North China Electric Power University\\
Yonghua Street 619, 071003, Baoding, P. R. China}
\email{tianjf\symbol{64}ncepu.edu.cn}

\begin{abstract}
In this paper, we investigate the recurrence relations for the Maclaurin coefficients of the products of elementary functions and hypergeometric functions. Specifically, we focus on the confluent hypergeometric function $\mathcal{M}(z) = h(z) M(a,c;z)$ and the Gaussian hypergeometric function $\mathcal{F}(z) = h(z) F(a,b;c;z)$, considering several specific choices for the function $h(z)$. In particular, we explore cases where $h(z)$ is chosen as $e^{pz}$, $(1-\theta z)^p$, $e^{-p \arctan z}$, $\sin(pz)$, $\cos(pz)$, $\sinh(pz)$, $\cosh(pz)$, $\arcsin(pz)$, and $\arccos(pz)$.
\end{abstract}

\thanks{*Corresponding author: Jing-Feng Tian(tianjf\symbol{64}ncepu.edu.cn)}

\maketitle

\tableofcontents

\section{Introduction}
A wide class of classical special functions arises as solutions to second-order linear ordinary differential equations of the form
\begin{equation*}
Y^{\prime\prime}(z) + P(z) Y^\prime(z) + Q(z) Y(z) = 0.
\end{equation*}
When, in a neighborhood of each singular point $z_0$, the quantities $(z-z_0) P(z)$ and $(z-z_0)^2 Q(z)$ remain bounded, the equation is referred to as a Fuchsian equation \cite{Yoshida-1987,Beyer-JCP-2021,Huang-JFA-2025}. A key feature of Fuchsian equations is that their solutions possess well-controlled behavior near singularities and admit local representations in terms of power series, possibly accompanied by a finite number of logarithmic factors.

The most classical and important examples include the confluent hypergeometric function (also known as the Kummer function) $M(a,c;z)$, and the Gaussian hypergeometric function $F(a,b;c;z)$, which are respectively defined as follows.
\begin{equation} \label{M-series}
M(a,c;z) = \sum_{n=0}^\infty \frac{(a)_n}{(c)_n \, n!} z^n, \quad -c \notin \mathbb{N} \cup \{0\}, a,c,z \in \mathbb{C},
\end{equation}
and
\begin{equation} \label{F-series}
F(a,b;c;z) = \sum_{n=0}^\infty \frac{(a)_n (b)_n}{(c)_n \, n!} z^n, \quad -c \notin \mathbb{N} \cup \{0\}, a,c,z \in \mathbb{C}, |z| <1,
\end{equation}
where $(a)_n = a(a+1)\cdots(a+n-1)$ is the Pochhammer symbol.

In the study of many special functions, one often encounters expressions in which a basic function is multiplied by a special function, of the form
\begin{equation*}
\mathcal{P}(z) := h(z) S(z),
\end{equation*}
where $h(z)$ denotes a basic function and $S(z)$ represents a special function.

%
%

\section{Recurrence relations for the Maclaurin coefficients of products of elementary functions and the confluent Hypergeometric functions}

In this section, we investigate the recurrence relations for the Maclaurin coefficients of functions of the form $\mathcal{P}(z) = h(z) M(a,c,z)$, where $h(z)$ is chosen as certain specific elementary functions.

\begin{thm} \label{thm-exp-M}
Let $a,c,p \in \mathbb{C}$ and $-c \notin \mathbb{N} \cup \{0\}$. Then
\begin{equation*}
e^{pz} M(a,c;z) = \sum_{n=0}^\infty u_n z^n, \quad z \in \mathbb{C},
\end{equation*}
if $u_0 = 1$, $u_1 = a/c+p$ and
\begin{equation*}
u_{n+1} = \frac{a+c p+2 n p+n}{(n+1) (c+n)} u_n -\frac{p (p+1)}{(n+1) (c+n)} u_{n-1}.
\end{equation*}
\end{thm}

Using the identities $\sinh(pz) = (e^{pz} - e^{-pz})/2$ and $\cosh(pz) = (e^{pz} + e^{-pz})/2$, Theorem~\ref{thm-exp-M} directly yields the following two theorems.

\begin{thm}
Let $a,c,p \in \mathbb{C}$ and $-c \notin \mathbb{N} \cup \{0\}$. Then
\begin{equation*}
\sinh(pz) M(a,c;z) = \sum_{n=0}^\infty \frac{u_n - v_n}{2} z^n, \quad z \in \mathbb{C},
\end{equation*}
if $u_0 = 1$, $u_1 = \frac{a}{c}+p$, $v_0 = 1$, $v_1 = \frac{a}{c}-p$,
\begin{equation*}
u_{n+1} = \frac{a+p(c+2 n)+n}{(n+1) (c+n)} u_n -\frac{p (p+1)}{(n+1) (c+n)} u_{n-1},
\end{equation*}
and
\begin{equation*}
v_{n+1} = \frac{a-p (c+2 n)+n}{(n+1) (c+n)} v_n -\frac{(p-1) p}{(n+1) (c+n)} v_{n-1},
\end{equation*}
\end{thm}

\begin{thm}
Let $a,c,p \in \mathbb{C}$ and $-c \notin \mathbb{N} \cup \{0\}$. Then
\begin{equation*}
\cosh(pz) M(a,c;z) = \sum_{n=0}^\infty \frac{u_n + v_n}{2} z^n, \quad z \in \mathbb{C},
\end{equation*}
if $u_0 = 1$, $u_1 = \frac{a}{c}+p$, $v_0 = 1$, $v_1 = \frac{a}{c}-p$,
\begin{equation*}
u_{n+1} = \frac{a+p(c+2 n)+n}{(n+1) (c+n)} u_n -\frac{p (p+1)}{(n+1) (c+n)} u_{n-1},
\end{equation*}
and
\begin{equation*}
v_{n+1} = \frac{a-p (c+2 n)+n}{(n+1) (c+n)} v_n -\frac{(p-1) p}{(n+1) (c+n)} v_{n-1},
\end{equation*}
\end{thm}

Observing that $\sin(pz) = (e^{ipz} - e^{-ipz})/(2i)$ and $\cos(pz) = (e^{ipz} + e^{-ipz})/2$, Theorem~\ref{thm-exp-M} directly yields the following two theorems.

\begin{thm}
Let $a,c,p \in \mathbb{C}$ and $-c \notin \mathbb{N} \cup \{0\}$. Then
\begin{equation*}
\sin(pz) M(a,c;z) = \sum_{n=0}^\infty \frac{u_n - v_n}{2i} z^n, \quad z \in \mathbb{C},
\end{equation*}
if $u_0 = 1$, $u_1 = \frac{a}{c}+ip$, $v_0 = 1$, $v_1 = \frac{a}{c}-ip$,
\begin{equation*}
u_{n+1} = \frac{a+ i p(c+2 n)+n}{(n+1) (c+n)} u_n -\frac{ip-p^2}{(n+1) (c+n)} u_{n-1},
\end{equation*}
and
\begin{equation*}
v_{n+1} = \frac{a-ip (c+2 n)+n}{(n+1) (c+n)} v_n + \frac{ip+p^2}{(n+1) (c+n)} v_{n-1},
\end{equation*}
\end{thm}

Next, we consider the product of $M(a,c;z)$ and $(1-\theta z)^p$, see also in Xu et al.\ \cite{Xu-Arxiv-2026}.

\begin{thm}
Let $a,c,\theta \in \mathbb{C}$, $-c \notin \mathbb{N} \cup \{0\}$, $p\in \mathbb{C}$. Then
\begin{equation*}
(1-\theta z)^p M(a,c;z) = \sum_{n=0}^\infty u_n z^n, \quad |z| < \frac{1}{|\theta|},
\end{equation*}
if $u_0 = 1$, $u_1 = \frac{a}{c}-\theta  p$, $u_2 = \frac{1}{2} \left(\frac{a^2+a}{c^2+c}-\frac{2 a \theta  p}{c}+\theta ^2 (p-1) p\right)$ and
\begin{equation*}
u_{n+1} = \beta_0(n) u_n + \beta_1(n) u_{n-1} + \beta_2(n) u_{n-2},
\end{equation*}
where
\begin{equation*}
\begin{aligned}
\beta_0(n) &= \frac{a+2 \theta  n (c+n-1)-\theta  p (c+2 n)+n}{(n+1) (c+n)},\\
\beta_1(n) &= \frac{\theta  (-2 a-\theta  (n-p-1) (c+n-p-2)-2 n+p+2)}{(n+1) (c+n)},\\
\beta_2(n) &= \frac{\theta ^2 (a+n-p-2)}{(n+1) (c+n)}.
\end{aligned}
\end{equation*}
\end{thm}

Next, we consider the product of $M(a,c;z)$ and $e^{-p \arctan z}$.

\begin{thm}
Let $a,c, p \in \mathbb{C}$ and $-c \notin \mathbb{N} \cup \{0\}$. Then
\begin{equation*}
e^{-p \arctan z} M(a,c;z) = \sum_{n=0}^\infty u_n z^n, \quad z \in \mathbb{C}.
\end{equation*}
if $u_0 = 1$, $u_1 = \frac{a}{c}- p$, $u_2 = \frac{1}{2} (\frac{a^2+a}{c^2+c}-\frac{2 a p}{c}+p^2)$, $u_3 = \frac{1}{6} (\frac{3 a p^2}{c}-\frac{3 a (a+1) p}{c (c+1)}+\frac{a (a+1) (a+2)}{c (c+1) (c+2)}-p^3+2 p)$,
\begin{equation*}
u_4= \frac{\left(
\begin{aligned}
& 6 a (a+1) (c+2) (c+3) p^2-4 a (c+1) (c+2) (c+3) \left(p^2-2\right) p\\
& \quad -4 a (a+1) (a+2) (c+3) p+a (a+1) (a+2) (a+3)\\
&\quad +c (c+1) (c+2) (c+3) \left(p^2-8\right) p^2
\end{aligned}\right)}{24 c (c+1) (c+2) (c+3)}
\end{equation*}
and
\begin{equation*}
u_{n+1} = \beta_0(n) u_n + \beta_1(n) u_{n-1} + \beta_2(n) u_{n-2} + \beta_3(n) u_{n-3} + \beta_4(n) u_{n-4},
\end{equation*}
where
\begin{equation*}
\begin{aligned}
\beta_0(n) &= \frac{a-p (c+2 n)+n}{(n+1) (c+n)},\\
\beta_1(n) &= \frac{-2 (n-1) (c+n-2)-p^2+p}{(n+1) (c+n)},\\
\beta_2(n) &= \frac{2 (a+n-2)-p (c+2 n-6)}{(n+1) (c+n)}, \\
\beta_3(n) &= \frac{p-(n-3) (c+n-4)}{(n+1) (c+n)}, \\
\beta_4(n) &= \frac{a+n-4}{(n+1) (c+n)}.
\end{aligned}
\end{equation*}
\end{thm}

Next, we revisit the product of trigonometric functions and the confluent hypergeometric function, and show that the corresponding coefficients satisfy a sixth-order recurrence relation.

\begin{thm} \label{thm-sin-M}
Let $a,c,p \in \mathbb{C}$ and $-c \notin \mathbb{N} \cup \{0\}$. Then
\begin{equation*}
\sin(pz) M(a,c;z) = \sum_{n=0}^\infty u_n z^n, \quad z \in \mathbb{C},
\end{equation*}
if $u_0=0$, $u_1=p$, $u_2=\frac{a p}{c}$, $u_3 = \frac{a (a+1) p}{2 c (c+1)}-\frac{p^3}{6}$, $u_4 = \frac{a (a+1) (a+2) p}{6 c (c+1) (c+2)}-\frac{a p^3}{6 c}$, $u_5=-\frac{a (a+1) p^3}{12 c (c+1)}+\frac{a (a+1) (a+2) (a+3) p}{24 c (c+1) (c+2) (c+3)}+\frac{p^5}{120}$,
\begin{equation*}
u_{n+1} = \sum_{i=0}^5 \beta_i(n) u_{n-i},
\end{equation*}
where
\begin{equation*}
\beta_0(n) = \frac{2 \left(a \left(c^2-2 c n+c-2 (n-2)^2\right)+c (n-1) (2 c+n-5)\right)}{(c-2) c (n+1) (c+n)},
\end{equation*}
\begin{equation*}
\beta_1(n) = \frac{\left(
\begin{aligned}
& -(n-4) (n-3) (n-2) (n-1) (4 p^2+1) + 2 (n-3) (n-2) (n-1) (4 a-c (4 p^2+3)) \\
& \quad +(n-2) (n-1) \left(8 a^2+a (4 c+6)-6 c ((c-2) p^2+c)+c\right) \\
& \quad + 2 (n-1) \left(a^2 (4 c-2)-3 a (c-1) c+c \left(-c^2+c+2\right) p^2\right) \\
& \quad -(c-2) \left(a^2 (c+2)-a c+c^2 (c+1) p^2\right)
\end{aligned}
\right)}{(c-2) c n (n+1) (c+n-1) (c+n)}, \\
\end{equation*}
\begin{equation*}
\beta_2(n) = \frac{\left(
\begin{aligned}
& -2 p^2 (c-3) \left(a (3 c+8)-2 c^2+c-32\right) \\
& \quad +2p^2 \left( n (-10 a+c (3 c-31)+104)+6 (c-6) n^2+4 n^3\right)\\
&\quad -2 (a+n-3) \left(2 a^2-a (c-2 n+3)-(n-2) (2 c+n-4)\right)
\end{aligned}
\right)}{(c-2) c n (n+1) (c+n-1) (c+n)}, \\
\end{equation*}
\begin{equation*}
\beta_3(n) = -\frac{\left(
\begin{aligned}
& p^2 \left(12 a^2-2 a (6 c+5)+c (6 c-1)\right)+2 (n-3) \left(a+c \left(4 p^2+3\right) p^2\right) \\
&\quad +(a-1) a+5 (c-2) c p^4+(n-4) (n-3) \left(8 p^4+6 p^2+1\right)
\end{aligned}
\right)}{(c-2) c n (n+1) (c+n-1) (c+n)}, \\
\end{equation*}
\begin{equation*}
\beta_4(n) = \frac{2 p^2 \left(p^2 (-6 a+5 c+4 (n-4))-3 a+2 c+n-4\right)}{(c-2) c n (n+1) (c+n-1) (c+n)},
\end{equation*}
\begin{equation*}
\beta_5(n) = -\frac{p^2 \left(4 p^4+5 p^2+1\right)}{(c-2) c n (n+1) (c+n-1) (c+n)}.
\end{equation*}
\end{thm}

\begin{thm}
Let $a,c,p \in \mathbb{C}$ and $-c \notin \mathbb{N} \cup \{0\}$. Then
\begin{equation*}
\cos(pz) M(a,c;z) = \sum_{n=0}^\infty u_n z^n, \quad z \in \mathbb{C},
\end{equation*}
if $u_0=1$, $u_1=\frac{a}{c}$, $u_2=\frac{1}{2} \left(\frac{a^2+a}{c^2+c}-p^2\right)$, $u_3 = \frac{a \left(\frac{(a+1) (a+2)}{(c+1) (c+2)}-3 p^2\right)}{6 c}$,
\begin{equation*}
u_4 = \frac{1}{24} \left(-\frac{6 a (a+1) p^2}{c (c+1)}+\frac{a (a+1) (a+2) (a+3)}{c (c+1) (c+2) (c+3)}+p^4\right),
\end{equation*}
\begin{equation*}
u_5=\frac{a \left(-\frac{10 (a+1) (a+2) p^2}{(c+1) (c+2)}+\frac{(a+1) (a+2) (a+3) (a+4)}{(c+1) (c+2) (c+3) (c+4)}+5 p^4\right)}{120 c},
\end{equation*}
and
\begin{equation*}
u_{n+1} = \sum_{i=0}^5 \beta_i(n) u_{n-i},
\end{equation*}
where $\beta_i(n)$ ($n=0,1,\cdots,5$) is defined in Theorem \ref{thm-sin-M}.
\end{thm}

Next, we revisit the product of hyperbolic functions and the confluent hypergeometric function, and show that the corresponding coefficients satisfy a sixth-order recurrence relation.

\begin{thm} \label{thm-sinh-M}
Let $a,c,p \in \mathbb{C}$ and $-c \notin \mathbb{N} \cup \{0\}$. Then
\begin{equation*}
\sinh(pz) M(a,c;z) = \sum_{n=0}^\infty u_n z^n, \quad z \in \mathbb{C},
\end{equation*}
if $u_0=0$,
$u_1=p$,
$u_2=\frac{a p}{c}$,
$u_3 = (a (1 + a) p)/(2 c (1 + c)) + p^3/6$,
$u_4 = \frac{a p^3}{6 c}+\frac{a (a+1) (a+2) p}{6 c (c+1) (c+2)}$,
$u_5= \frac{a (a+1) p^3}{12 c (c+1)}+\frac{a (a+1) (a+2) (a+3) p}{24 c (c+1) (c+2) (c+3)}+\frac{p^5}{120}$,
\begin{equation*}
u_{n+1} = \sum_{i=0}^5 \beta_i(n) u_{n-i},
\end{equation*}
where
\begin{equation*}
\beta_0(n) = \frac{2 \left(a \left(c^2-2 c n+c-2 (n-2)^2\right)+c (n-1) (2 c+n-5)\right)}{(c-2) c (n+1) (c+n)},
\end{equation*}
\begin{equation*}
\beta_1(n) = \frac{\left(
\begin{aligned}
& (n-4) (n-3) (n-2) (n-1) (4 p^2-1) \\
& \quad +2 (n-3) (n-2) (n-1) \left(4 a+c \left(4 p^2-3\right)\right) \\
& \quad +(n-2) (n-1) \left(8 a^2+a (4 c+6)+c \left(6 (c-2) p^2-6 c+1\right)\right) \\
& \quad +2 (n-1) \left(a^2 (4 c-2)-3 a (c-1) c+(c-2) c (c+1) p^2\right) \\
& \quad +(c-2) \left(a^2 (-(c+2))+a c+c^2 (c+1) p^2\right)
\end{aligned}
\right)}{(c-2) c n (n+1) (c+n-1) (c+n)}, \\
\end{equation*}
\begin{equation*}
\beta_2(n) = \frac{\left(
\begin{aligned}
& 2 p^2 \left((c-3) \left(a (3 c+8)-2 c^2+c-32\right) \right)\\
&\quad  +2 p^2 \left(n (10 a+c (31-3 c)-104)-6 (c-6) n^2-4 n^3\right) \\
& \quad -2 (a+n-3) \left(2 a^2-a (c-2 n+3)-(n-2) (2 c+n-4)\right)
\end{aligned}
\right)}{(c-2) c n (n+1) (c+n-1) (c+n)}, \\
\end{equation*}
\begin{equation*}
\beta_3(n) = -\frac{\left(
\begin{aligned}
& p^2 \left(-12 a^2+2 a (6 c+5)-6 c^2+c\right)+2 (n-3) \left(a+c \left(4 p^2-3\right) p^2\right) \\
& \quad +(a-1) a+5 (c-2) c p^4+(n-4) (n-3) \left(8 p^4-6 p^2+1\right)
\end{aligned}
\right)}{(c-2) c n (n+1) (c+n-1) (c+n)}, \\
\end{equation*}
\begin{equation*}
\beta_4(n) = \frac{2 p^2 \left(p^2 (-6 a+5 c+4 (n-4))+3 a-2 c-n+4\right)}{(c-2) c n (n+1) (c+n-1) (c+n)},
\end{equation*}
\begin{equation*}
\beta_5(n) = \frac{4 p^6-5 p^4+p^2}{(c-2) c n (n+1) (c+n-1) (c+n)}.
\end{equation*}
\end{thm}

\begin{thm}
Let $a,c,p \in \mathbb{C}$ and $-c \notin \mathbb{N} \cup \{0\}$. Then
\begin{equation*}
\cosh(pz) M(a,c;z) = \sum_{n=0}^\infty u_n z^n, \quad z \in \mathbb{C},
\end{equation*}
if $u_0=1$,
$u_1=\frac{a}{c}$,
$u_2=\frac{a (a+1)}{2 c (c+1)}+\frac{p^2}{2}$,
$u_3 = \frac{a p^2}{2 c}+\frac{a (a+1) (a+2)}{6 c (c+1) (c+2)}$,
$u_4 = \frac{a (a+1) p^2}{4 c (c+1)}+\frac{a (a+1) (a+2) (a+3)}{24 c (c+1) (c+2) (c+3)}+\frac{p^4}{24}$,
$u_5= \frac{a p^4}{24 c}+\frac{a (a+1) (a+2) p^2}{12 c (c+1) (c+2)}+\frac{a (a+1) (a+2) (a+3) (a+4)}{120 c (c+1) (c+2) (c+3) (c+4)}$,
\begin{equation*}
u_{n+1} = \sum_{i=0}^5 \beta_i(n) u_{n-i},
\end{equation*}
where $\beta_i(n)$ ($n=0,1,\cdots,5$) is defined in Theorem \ref{thm-sinh-M}.
\end{thm}

At the end of this section, we investigate the product of inverse trigonometric functions and the confluent hypergeometric function.

\begin{thm} \label{thm-arcsin-M}
Let $a,c,p \in \mathbb{C}$ and $-c \notin \mathbb{N} \cup \{0\}$. Then
\begin{equation*}
\arcsin(pz) M(a,c;z) = \sum_{n=0}^\infty u_n z^n, \quad z \in \mathbb{C},
\end{equation*}
if $u_0=0$, $u_1=p$, $u_2=\frac{a p}{c}$, $u_3=\frac{a (a+1) p}{2 c (c+1)}+\frac{p^3}{6}$, $u_4=\frac{a p^3}{6 c}+\frac{a (a+1) (a+2) p}{6 c (c+1) (c+2)}$,
\begin{equation*}
u_5=\frac{a (a+1) p^3}{12 c (c+1)}+\frac{a (a+1) (a+2) (a+3) p}{24 c (c+1) (c+2) (c+3)}+\frac{3 p^5}{40},
\end{equation*}
\begin{equation*}
u_6=\frac{3 a p^5}{40 c}+\frac{a (a+1) (a+2) p^3}{36 c (c+1) (c+2)}+\frac{a (a+1) (a+2) (a+3) (a+4) p}{120 c (c+1) (c+2) (c+3) (c+4)},
\end{equation*}
\begin{equation*}
u_7=\frac{3 a (a+1) p^5}{80 c (c+1)}+\frac{a (a+1) (a+2) (a+3) p^3}{144 c (c+1) (c+2) (c+3)}+\frac{a (a+1) (a+2) (a+3) (a+4) (a+5) p}{720 c (c+1) (c+2) (c+3) (c+4) (c+5)}+\frac{5 p^7}{112},
\end{equation*}
\begin{equation*}
\begin{aligned}
u_8 & =\frac{5 a p^7}{112 c}+\frac{a (a+1) (a+2) p^5}{80 c (c+1) (c+2)}+\frac{a (a+1) (a+2) (a+3) (a+4) p^3}{720 c (c+1) (c+2) (c+3) (c+4)}\\
&\quad +\frac{a (a+1) (a+2) (a+3) (a+4) (a+5) (a+6) p}{5040 c (c+1) (c+2) (c+3) (c+4) (c+5) (c+6)},
\end{aligned}
\end{equation*}
\begin{equation*}
\begin{aligned}
u_9&=\frac{5 a (a+1) p^7}{224 c (c+1)}+\frac{a (a+1) (a+2) (a+3) p^5}{320 c (c+1) (c+2) (c+3)} +\frac{a (a+1) (a+2) (a+3) (a+4) (a+5) p^3}{4320 c (c+1) (c+2) (c+3) (c+4) (c+5)}\\
&\quad +\frac{a (a+1) (a+2) (a+3) (a+4) (a+5) (a+6) (a+7) p}{40320 c (c+1) (c+2) (c+3) (c+4) (c+5) (c+6) (c+7)}+\frac{35 p^9}{1152},
\end{aligned}
\end{equation*}
\begin{equation*}
\begin{aligned}
u_{10}& =\frac{35 a p^9}{1152 c}+\frac{5 a (a+1) (a+2) p^7}{672 c (c+1) (c+2)}+\frac{a (a+1) (a+2) (a+3) (a+4) p^5}{1600 c (c+1) (c+2) (c+3) (c+4)} \\
& \quad +\frac{a (a+1) (a+2) (a+3) (a+4) (a+5) (a+6) p^3}{30240 c (c+1) (c+2) (c+3) (c+4) (c+5) (c+6)} \\
& \quad +\frac{a (a+1) (a+2) (a+3) (a+4) (a+5) (a+6) (a+7) (a+8) p}{362880 c (c+1) (c+2) (c+3) (c+4) (c+5) (c+6) (c+7) (c+8)},
\end{aligned}
\end{equation*}
\begin{equation*}
\begin{aligned}
u_{11}&=\frac{35 a (a+1) p^9}{2304 c (c+1)}+\frac{5 a (a+1) (a+2) (a+3) p^7}{2688 c (c+1) (c+2) (c+3)}+\frac{a (a+1) (a+2) (a+3) (a+4) (a+5) p^5}{9600 c (c+1) (c+2) (c+3) (c+4) (c+5)} \\
& \quad +\frac{a (a+1) (a+2) (a+3) (a+4) (a+5) (a+6) (a+7) p^3}{241920 c (c+1) (c+2) (c+3) (c+4) (c+5) (c+6) (c+7)} \\
& \quad +\frac{a (a+1) (a+2) (a+3) (a+4) (a+5) (a+6) (a+7) (a+8) (a+9) p}{3628800 c (c+1) (c+2) (c+3) (c+4) (c+5) (c+6) (c+7) (c+8) (c+9)}+\frac{63 p^{11}}{2816},
\end{aligned}
\end{equation*}
and
\begin{equation*}
u_{n+1} = \sum_{i=0}^{11} \beta_i(n) u_{n-i}, n\geq11,
\end{equation*}
where
\begin{equation*}
\beta_0(n) = \frac{2 a \left(c^2-2 c n+c-2 (n-2)^2\right)+(c-2) c p^2 (c+2 n-1)+2 c (n-1) (2 c+n-5)}{(c-2) c (n+1) (c+n)},
\end{equation*}
\begin{equation*}
\beta_1(n) = \frac{\left(
\begin{aligned}
& (n-4) (n-3) (n-2) (n-1) \left(4 (c-2) c p^2+p^4-1\right) \\
& \quad +2 (n-3) (n-2) (n-1) \left(4 a \left(p^2+1\right)+c \left(2 (2 (c-1) c-5) p^2+p^4-3\right)\right) \\
& \quad +(n-2) (n-1) \left(8 a^2+p^2 (6 a (2 c+1)+c (c (4 (c-1) c-17)+3))+4 a c+6 a-6 c^2+c\right) \\
& \quad +(n-1) \left(a^2 (8 c-4)+2 a c \left(c \left(p^2-3\right)+p^2+3\right)-c (c+1) p^2 \left(c \left(p^2+4\right)-2 \left(p^2+3\right)\right)\right) \\
& \quad -a (c-2) \left(a (c+2)+c \left((c+1) p^2-1\right)\right)
\end{aligned}\right)}{(c-2) c n (n+1) (c+n-1) (c+n)},
\end{equation*}
\begin{equation*}
\beta_2(n) = -\frac{\left(
\begin{aligned}
& 4 (n-5) (n-4) (n-3) (n-2) p^2 \left(-4 a+2 c+p^2\right) \\
& \quad +2 (n-4) (n-3) (n-2) \left(
\begin{aligned}
& -p^2 (8 a (2 c+1)+4 c (1-3 c)+1) \\
& \quad +(c-1) (3 c+1) p^4+p^6-1
\end{aligned}\right) \\
& \quad -(n-3) (n-2) \left(
\begin{aligned}
& c \left((24-9 (c-1) c) p^4+(8 c (1-2 c)+33) p^2-3 p^6+4\right) \\
& \quad -2 a \left((6-4 c (c+1)) p^2+3 p^4+1\right)
\end{aligned}\right) \\
& \quad +(n-2) \left(
\begin{aligned}
& p^2 \left(8 a^2+4 a c (2 (c-1) c-3)+6 a-6 c^2+c\right) \\
& \quad +p^4 (a (6 c-2)+c (c (c (3 c-7)-7)+13))\\
& \quad +2 a (4 a-3 c+2)+(c-2) c p^6
\end{aligned}\right) \\
& \quad +a \left(p^2 (a (4 c-2)-3 (c-1) c)+2 (a-1) (2 a-c+1)-(c-2) (c+1) p^4\right)
\end{aligned}\right)}{(c-2) c n (n+1) (c+n-1) (c+n)},
\end{equation*}
\begin{equation*}
\beta_3(n) = -\frac{\left(
\begin{aligned}
& 2 (n-6) (n-5) (n-4) (n-3) p^2 \left(3 (c-2) c p^2+p^4-2\right) \\
&\quad +4 (n-5) (n-4) (n-3) p^2 \left(p^2 (6 a+3 c ((c-1) c-3)-2)+8 a+(c-1) p^4-6 c\right) \\
& \quad +(n-4) (n-3) \left(
\begin{aligned}
& p^2 \left(32 a^2+8 a (2 c+3)+4 c (1-6 c)+3\right) \\
&\quad +p^4 (6 a (6 c-1)+3 (c-3) c (2 c (c+2)-1)+6)\\
&\quad +(3-6 c) p^6-p^8+1
\end{aligned}\right) \\
& \quad +(n-3) \left(
\begin{aligned}
& -2 p^2 \left(a^2 (8-16 c)+12 a (c-1) c+a-2 c\right)-2 p^6 \left(a+c^3-6 c\right) \\
&\quad +3 p^4 \left(2 a ((c-3) c+1)+c \left(-4 c^2+6 c+7\right)\right)+2 a-c p^8
\end{aligned}\right) \\
& \quad +a p^2 \left(-4 a \left(c^2-3\right)+c (4 c-5)-2\right)\\
&\quad + a \left(p^4 (-5 a+c (c (7-3 c)+6)-11)+a-(c-2) p^6-1\right)
\end{aligned}\right)}{(c-2) c n (n+1) (c+n-1) (c+n)},
\end{equation*}
\begin{equation*}
\beta_4(n) = \frac{p^2\left(
\begin{aligned}
& 4 (n-7) (n-6) (n-5) (n-4) p^2 \left(-6 a+3 c+2 p^2\right) \\
& \quad +2 (n-6) (n-5) (n-4) \left(
\begin{aligned}
& -3 p^2 \left(a (8 c+4)-6 c^2+2 c+1\right) \\
&\quad +(c (3 c+2)+6) p^4+p^6-4
\end{aligned}\right) \\
& \quad +(n-5) (n-4) \left(
\begin{aligned}
& -3 p^2 \left(4 a \left(c^2+c-3\right)+c \left(-8 c^2+4 c+21\right)\right) \\
&\quad +3 p^4 (4 a+c (3 (c-1) c-2)-2)+8 (a-2 c)+(3 c-2) p^6
\end{aligned}\right) \\
& \quad +(n-4) \left(
\begin{aligned}
& 3 p^2 \left(8 a^2+a (4 c ((c-1) c-1)-2)+c (5-6 c)+1\right)+32 a^2 \\
&\quad +p^4 \left(4 a (3 c-2)+(3 c-11) c^3+6 c+5\right)\\
&\quad -24 a c+16 a+((c-4) c-1) p^6+1
\end{aligned}\right) \\
& \quad -a \left(-16 a^2+p^2 (a (26-12 c)+c (9 c-17)-7) \right) \\
&\quad -a \left(8 a c+8 a+\left(2 c^2-13\right) p^4-8 c+p^6+7\right)
\end{aligned}\right)}{(c-2) c n (n+1) (c+n-1) (c+n)},
\end{equation*}
\begin{equation*}
\beta_5(n) = \frac{\left(
\begin{aligned}
& (n-8) (n-7) (n-6) (n-5) p^4 \left(4 (c-2) c p^2+p^4-6\right) \\
&\quad +2 (n-7) (n-6) (n-5) p^4 \left(2 p^2 (6 a+c (2 (c-1) c-7)-4)+6 (4 a-3 c)+(c-2) p^4\right) \\
&\quad -(n-6) (n-5) p^2 \left(
\begin{aligned}
& -3 p^2 \left(8 a c+4 a (4 a+3)-12 c^2+2 c+3\right) \\
&\quad +p^4 (6 a (5-6 c)+c (c (35-4 (c-1) c)-9)+12) \\
& \quad +(6 c+5) p^6-4
\end{aligned}\right) \\
&\quad +(n-5) \left(
\begin{aligned}
& p^8 \left(-\left(2 a+(c+1) \left(c^2-2\right)\right)\right) \\
&\quad +2 p^6 (a (3 (c-7) c-2)+c (3 c (5-2 c)+4)+1)\\
&\quad +2 p^4 (3 a (a (8 c-4)-6 (c-1) c-1)+6 c-2)+8 a p^2
\end{aligned}\right) \\
&\quad -a p^2 \left(
\begin{aligned}
& p^2 \left(6 (a-1) c^2-12 a+3 c+10\right) \\
&\quad +p^4 (10 a+c (c (3 c-11)+2)+4)-4 a+(c-4) p^6+4
\end{aligned}\right)
\end{aligned}\right)}{(c-2) c n (n+1) (c+n-1) (c+n)},
\end{equation*}
\begin{equation*}
\beta_6(n) = \frac{p^4\left(
\begin{aligned}
& -4 (n-9) (n-8) (n-7) (n-6) p^2 \left(-4 a+2 c+p^2\right) \\
&\quad -2 (n-8) (n-7) (n-6) \left(-p^2 (8 a (2 c+1)+4 c (1-3 c)+3)+(c (c+2)+7) p^4-6\right) \\
&\quad -(n-7) (n-6) \left(
\begin{aligned}
& 3 p^4 (2 a+c ((c-1) c+4)-2) \\
&\quad -p^2 (4 a (2 c (c+1)-9)+c (8 c (1-2 c)+51))+12 (a-2 c)
\end{aligned}\right) \\
&\quad -(n-6) \left(
\begin{aligned}
& p^2 \left(24 a^2+a (4 c (2 (c-1) c-1)-30)+9 c (3-2 c)-2\right) \\
&\quad +p^4 (6 a (c-1)+c (c ((c-5) c+5)-7)-3)+12 a (4 a-3 c+2)+3
\end{aligned}\right) \\
&\quad +a \left(
\begin{aligned}
& p^2 (a (46-12 c)+(c-3) (9 c+2)) \\
&\quad +3 (4 a (-2 a+c+1)-4 c+3)+\left(c^2+c-3\right) p^4
\end{aligned}\right)
\end{aligned}\right)}{(c-2) c n (n+1) (c+n-1) (c+n)},
\end{equation*}
\begin{equation*}
\beta_7(n) = -\frac{\left(
\begin{aligned}
& (n-10) (n-9) (n-8) (n-7) p^6 \left((c-2) c p^2-4\right) \\
&\quad +2 (n-9) (n-8) (n-7) p^6 \left(4 a \left(p^2+4\right)+(c ((c-1) c-4)-4) p^2-12 c\right) \\
&\quad +(n-8) (n-7) p^4 \left(
\begin{aligned}
& p^2 \left(32 a^2+8 a (2 c+3)+4 c (1-6 c)+9\right) \\
&\quad +p^4 (6 a (2 c-3)+c (c+3) ((c-4) c+1)-18)+6
\end{aligned}\right) \\
&\quad +(n-7) \left(
\begin{aligned}
& p^8 (2 a ((c-11) c-5)+c (2 (7-2 c) c-7)+2) \\
&\quad +2 p^6 (a (2 c-1) (8 a-6 c+3)+6 c-4)+12 a p^4
\end{aligned}\right) \\
&\quad +a p^4 \left(
\begin{aligned}
& p^2 \left(-4 (a-1) c^2+4 a+c\right)+p^4 (-(5 a+(c-3) (c-2) c-7)) \\
&\quad +6 a-2 \left(7 p^2+3\right)
\end{aligned}\right)
\end{aligned}\right)}{(c-2) c n (n+1) (c+n-1) (c+n)},
\end{equation*}
\begin{equation*}
\beta_8(n) = \frac{p^6\left(
\begin{aligned}
& -2 (n-11) (n-10) (n-9) (n-8) p^2 (2 a-c) \\
&\quad -2 (n-10) (n-9) (n-8) \left(p^2 \left(a (4 c+2)-3 c^2+c+1\right)+4\right) \\
&\quad -(n-9) (n-8) \left(
p^2 \left(2 a \left(c^2+c-6\right) +c \left(-4 c^2+2 c+15\right)\right) -8 a+16 c
\right) \\
&\quad +(n-8) \left(
\begin{aligned}
& 32 a^2+p^2 \left(2 a \left(4 a+(c-1) c^2-9\right) +c (13-6 c)-5\right)\\
&\quad -24 a c+16 a+3
\end{aligned}\right) \\
&\quad -a \left(p^2 (a (22-4 c)+c (3 c-11))+8 a (-2 a+c+1)-8 c+p^2+5\right)
\end{aligned}\right)}{(c-2) c n (n+1) (c+n-1) (c+n)},
\end{equation*}
\begin{equation*}
\beta_9(n) = -\frac{p^6\left(
\begin{aligned}
& p^2 \left(
\begin{aligned}
& a \left(c^2 (6 n-55)+c \left(-4 n^2+70 n-307\right)-8 n^3+234 n^2-2276 n+7368\right) \\
& \quad +a^2 \left(c^2-8 c (n-9)-8 n^2+156 n-756\right) \\
&\quad +(n-9) \left(6 c^2 (n-10)+c \left(6 n^2-127 n+666\right)+n^3-33 n^2+359 n-1286\right)
\end{aligned}\right) \\
& \quad - 4 (a+n-10) (a+n-9) \\
\end{aligned}\right)}{(c-2) c n (n+1) (c+n-1) (c+n)},
\end{equation*}
\begin{equation*}
\beta_{10}(n) = \frac{p^8 \left(
\begin{aligned}
& -4 a^3+2 a^2 (c-4 n+41)+a \left(c (6 n-62)-2 n^2+38 n-179\right) \\
&\quad +(n-10) \left(4 c (n-11)+2 n^2-46 n+263\right)
\end{aligned}\right)}{(c-2) c n (n+1) (c+n-1) (c+n)},
\end{equation*}
\begin{equation*}
\beta_{11}(n) = -\frac{p^8 (a+n-12) (a+n-11)}{(c-2) c n (n+1) (c+n-1) (c+n)}.
\end{equation*}
\end{thm}

\begin{thm} \label{thm-arccos-M}
Let $a,c,p \in \mathbb{C}$ and $-c \notin \mathbb{N} \cup \{0\}$. Then
\begin{equation*}
\arccos(pz) M(a,c;z) = \sum_{n=0}^\infty u_n z^n, \quad z \in \mathbb{C},
\end{equation*}
if $u_0=\frac{\pi }{2}$, $u_1=\frac{\pi  a}{2 c}-p$, $u_2=\frac{\pi  a (a+1)}{4 c (c+1)}-\frac{a p}{c}$, $u_3=-\frac{a (a+1) p}{2 c (c+1)}+\frac{\pi  a (a+1) (a+2)}{12 c (c+1) (c+2)}-\frac{p^3}{6}$,
\begin{equation*}
u_4=-\frac{a p^3}{6 c}-\frac{a (a+1) (a+2) p}{6 c (c+1) (c+2)}+\frac{\pi  a (a+1) (a+2) (a+3)}{48 c (c+1) (c+2) (c+3)},
\end{equation*}
\begin{equation*}
u_5=-\frac{a (a+1) p^3}{12 c (c+1)}-\frac{a (a+1) (a+2) (a+3) p}{24 c (c+1) (c+2) (c+3)}+\frac{\pi  a (a+1) (a+2) (a+3) (a+4)}{240 c (c+1) (c+2) (c+3) (c+4)}-\frac{3 p^5}{40},
\end{equation*}
\begin{equation*}
\begin{aligned}
u_6 &= -\frac{3 a p^5}{40 c}-\frac{a (a+1) (a+2) p^3}{36 c (c+1) (c+2)}-\frac{a (a+1) (a+2) (a+3) (a+4) p}{120 c (c+1) (c+2) (c+3) (c+4)}\\
& \quad +\frac{\pi  a (a+1) (a+2) (a+3) (a+4) (a+5)}{1440 c (c+1) (c+2) (c+3) (c+4) (c+5)},
\end{aligned}
\end{equation*}
\begin{equation*}
\begin{aligned}
u_7&=-\frac{3 a (a+1) p^5}{80 c (c+1)}-\frac{a (a+1) (a+2) (a+3) p^3}{144 c (c+1) (c+2) (c+3)}-\frac{a (a+1) (a+2) (a+3) (a+4) (a+5) p}{720 c (c+1) (c+2) (c+3) (c+4) (c+5)} \\
& \quad +\frac{\pi  a (a+1) (a+2) (a+3) (a+4) (a+5) (a+6)}{10080 c (c+1) (c+2) (c+3) (c+4) (c+5) (c+6)}-\frac{5 p^7}{112},
\end{aligned}
\end{equation*}
\begin{equation*}
\begin{aligned}
u_8 & = -\frac{5 a p^7}{112 c}-\frac{a (a+1) (a+2) p^5}{80 c (c+1) (c+2)}-\frac{a (a+1) (a+2) (a+3) (a+4) p^3}{720 c (c+1) (c+2) (c+3) (c+4)}\\
& \quad -\frac{a (a+1) (a+2) (a+3) (a+4) (a+5) (a+6) p}{5040 c (c+1) (c+2) (c+3) (c+4) (c+5) (c+6)} \\
& \quad +\frac{\pi  a (a+1) (a+2) (a+3) (a+4) (a+5) (a+6) (a+7)}{80640 c (c+1) (c+2) (c+3) (c+4) (c+5) (c+6) (c+7)},
\end{aligned}
\end{equation*}
\begin{equation*}
\begin{aligned}
u_9 &=-\frac{5 a (a+1) p^7}{224 c (c+1)}-\frac{a (a+1) (a+2) (a+3) p^5}{320 c (c+1) (c+2) (c+3)}-\frac{a (a+1) (a+2) (a+3) (a+4) (a+5) p^3}{4320 c (c+1) (c+2) (c+3) (c+4) (c+5)} \\
& \quad -\frac{a (a+1) (a+2) (a+3) (a+4) (a+5) (a+6) (a+7) p}{40320 c (c+1) (c+2) (c+3) (c+4) (c+5) (c+6) (c+7)} \\
& \quad +\frac{\pi  a (a+1) (a+2) (a+3) (a+4) (a+5) (a+6) (a+7) (a+8)}{725760 c (c+1) (c+2) (c+3) (c+4) (c+5) (c+6) (c+7) (c+8)}-\frac{35 p^9}{1152},
\end{aligned}
\end{equation*}
\begin{equation*}
\begin{aligned}
u_{10}& =-\frac{35 a p^9}{1152 c}-\frac{5 a (a+1) (a+2) p^7}{672 c (c+1) (c+2)}-\frac{a (a+1) (a+2) (a+3) (a+4) p^5}{1600 c (c+1) (c+2) (c+3) (c+4)} \\
& \quad -\frac{a (a+1) (a+2) (a+3) (a+4) (a+5) (a+6) p^3}{30240 c (c+1) (c+2) (c+3) (c+4) (c+5) (c+6)}\\
& \quad-\frac{a (a+1) (a+2) (a+3) (a+4) (a+5) (a+6) (a+7) (a+8) p}{362880 c (c+1) (c+2) (c+3) (c+4) (c+5) (c+6) (c+7) (c+8)}\\
& \quad +\frac{\pi  a (a+1) (a+2) (a+3) (a+4) (a+5) (a+6) (a+7) (a+8) (a+9)}{7257600 c (c+1) (c+2) (c+3) (c+4) (c+5) (c+6) (c+7) (c+8) (c+9)},
\end{aligned}
\end{equation*}
\begin{equation*}
\begin{aligned}
u_{11}&=-\frac{35 a (a+1) p^9}{2304 c (c+1)}-\frac{5 a (a+1) (a+2) (a+3) p^7}{2688 c (c+1) (c+2) (c+3)}-\frac{a (a+1) (a+2) (a+3) (a+4) (a+5) p^5}{9600 c (c+1) (c+2) (c+3) (c+4) (c+5)}\\
& \quad -\frac{a (a+1) (a+2) (a+3) (a+4) (a+5) (a+6) (a+7) p^3}{241920 c (c+1) (c+2) (c+3) (c+4) (c+5) (c+6) (c+7)}\\
& \quad-\frac{a (a+1) (a+2) (a+3) (a+4) (a+5) (a+6) (a+7) (a+8) (a+9) p}{3628800 c (c+1) (c+2) (c+3) (c+4) (c+5) (c+6) (c+7) (c+8) (c+9)}\\
& \quad+\frac{\pi  a (a+1) (a+2) (a+3) (a+4) (a+5) (a+6) (a+7) (a+8) (a+9) (a+10)}{79833600 c (c+1) (c+2) (c+3) (c+4) (c+5) (c+6) (c+7) (c+8) (c+9) (c+10)}-\frac{63 p^{11}}{2816},
\end{aligned}
\end{equation*}
and
\begin{equation*}
u_{n+1} = \sum_{i=0}^{11} \beta_i(n) u_{n-i}, n\geq11,
\end{equation*}
where $\beta_i(n)$ ($n=0,1,\cdots,11$) is defined in Theorem \ref{thm-arcsin-M}.
\end{thm}

\section{Recurrence relations for the Maclaurin coefficients of products of elementary functions and the Gauss Hypergeometric functions}

We first consider the product of the exponential function and the Gauss hypergeometric function.

\begin{thm} \label{thm-exp-F}
Let $a,b,c,p,z \in \mathbb{C}$ and $-c \notin \mathbb{N} \cup \{0\}$. Suppose that
\begin{equation*}
e^{pz} F(a,b;c;z) = \sum_{n=0}^\infty u_n z^n, \quad |z| < 1.
\end{equation*}
Then $u_0 = 1$, $u_1 = \frac{ab}{c}+p$, $u_2=((a (1 + a) b (1 + b))/(2 c (1 + c)) + (a b p)/c + p^2/2)$ and
\begin{equation*}
u_{n+1} = \frac{(a+n) (b+n)+p (c+2 n)}{(n+1) (c+n)} u_n -\frac{p (a+b+2 n+p-1)}{(n+1) (c+n)} u_{n-1} + \frac{p^2}{(n+1) (c+n)} u_{n-2}.
\end{equation*}
\end{thm}

\begin{remark}
Jian Shi from Hebei University mentioned, in a WeChat communication with the corresponding author of this paper, that he has derived recurrence relations for the Maclaurin series of the product of exponential and Gaussian hypergeometric functions over the real domain. These results have not yet been made available to us, and we have not been able to find such results in the publicly available literature.
\end{remark}

Noting that $\sinh(pz) = (e^{pz} - e^{-pz})/2$ and $\cosh(pz) = (e^{pz} + e^{-pz})/2$, we have
\begin{equation*}
\sinh(pz) F(a,b;c;z) = \frac{1}{2} \Big( e^{pz} F(a,b;c;z) - e^{-pz} F(a,b;c;z) \Big),
\end{equation*}
and
\begin{equation*}
\cosh(pz) F(a,b;c;z) = \frac{1}{2} \Big( e^{pz} F(a,b;c;z) + e^{-pz} F(a,b;c;z) \Big),
\end{equation*}
which leads to the following two theorems from Theorem \ref{thm-exp-F}.

\begin{thm} \label{thm-sinh-F-2}
Let $a,b,c,p,z \in \mathbb{C}$ and $-c \notin \mathbb{N} \cup \{0\}$. Suppose that
\begin{equation*}
\sinh(pz) F(a,b;c;z) = \sum_{n=0}^\infty \frac{u_n - v_n}{2} z^n, \quad |z|<1.
\end{equation*}
Then $u_0 = 1$, $u_1 = \frac{ab}{c}+p$, $u_2=((a (1 + a) b (1 + b))/(2 c (1 + c)) + (a b p)/c + p^2/2)$, $v_0 = 1$, $v_1 = \frac{ab}{c}-p$, $u_2=((a (1 + a) b (1 + b))/(2 c (1 + c)) - (a b p)/c + p^2/2)$
\begin{equation*}
u_{n+1} = \frac{(a+n) (b+n)+p (c+2 n)}{(n+1) (c+n)} u_n -\frac{p (a+b+2 n+p-1)}{(n+1) (c+n)} u_{n-1} + \frac{p^2}{(n+1) (c+n)} u_{n-2},
\end{equation*}
and
\begin{equation*}
v_{n+1} = \frac{(a+n) (b+n)-p (c+2 n)}{(n+1) (c+n)} v_n + \frac{p (a+b+2 n-p-1)}{(n+1) (c+n)} v_{n-1} + \frac{p^2}{(n+1) (c+n)} v_{n-2}.
\end{equation*}
\end{thm}

\begin{thm} \label{thm-cosh-F-2}
Let $a,b,c,p,z \in \mathbb{C}$ and $-c \notin \mathbb{N} \cup \{0\}$. Suppose that
\begin{equation*}
\cosh(pz) F(a,b;c;z) = \sum_{n=0}^\infty \frac{u_n + v_n}{2} z^n, \quad |z|<1.
\end{equation*}
Then $u_0 = 1$, $u_1 = \frac{ab}{c}+p$, $u_2=((a (1 + a) b (1 + b))/(2 c (1 + c)) + (a b p)/c + p^2/2)$, $v_0 = 1$, $v_1 = \frac{ab}{c}-p$, $u_2=((a (1 + a) b (1 + b))/(2 c (1 + c)) - (a b p)/c + p^2/2)$
\begin{equation*}
u_{n+1} = \frac{(a+n) (b+n)+p (c+2 n)}{(n+1) (c+n)} u_n -\frac{p (a+b+2 n+p-1)}{(n+1) (c+n)} u_{n-1} + \frac{p^2}{(n+1) (c+n)} u_{n-2},
\end{equation*}
and
\begin{equation*}
v_{n+1} = \frac{(a+n) (b+n)-p (c+2 n)}{(n+1) (c+n)} v_n + \frac{p (a+b+2 n-p-1)}{(n+1) (c+n)} v_{n-1} + \frac{p^2}{(n+1) (c+n)} v_{n-2}.
\end{equation*}
\end{thm}

Noting that $\sin(pz) = \frac{e^{ipz} - e^{-ipz}}{2i}$ and $\cos(pz) = \frac{e^{ipz} + e^{-ipz}}{2}$, we have
\begin{equation*}
\sin(pz) F(a,b;c;z) = \frac{1}{2i} \Big( e^{ipz} F(a,b;c;z) - e^{-ipz} F(a,b;c;z) \Big),
\end{equation*}
and
\begin{equation*}
\cos(pz) F(a,b;c;z) = \frac{1}{2} \Big( e^{ipz} F(a,b;c;z) + e^{-ipz} F(a,b;c;z) \Big),
\end{equation*}
which leads to the following two theorems from Theorem \ref{thm-exp-F}.

\begin{thm} \label{thm-sin-F-2}
Let $a,b,c,p,z \in \mathbb{C}$ and $-c \notin \mathbb{N} \cup \{0\}$. Suppose that
\begin{equation*}
\sin(pz) F(a,b;c;z) = \sum_{n=0}^\infty \frac{u_n - v_n}{2 i} z^n, \quad |z|<1.
\end{equation*}
Then $u_0 = 1$, $u_1 = \frac{ab}{c}+ip$, $u_2=\frac{i a b p}{c}+\frac{a (a+1) b (b+1)}{2 c (c+1)}-\frac{p^2}{2}$, $v_0 = 1$, $v_1 = \frac{ab}{c}-ip$, $u_2=-\frac{i a b p}{c}+\frac{a (a+1) b (b+1)}{2 c (c+1)}-\frac{p^2}{2}$
\begin{equation*}
u_{n+1} = \frac{(a+n) (b+n)+i p (c+2 n)}{(n+1) (c+n)} u_n + \frac{p (p-i (a+b+2 n-1))}{(n+1) (c+n)} u_{n-1} -\frac{p^2}{(n+1) (c+n)} u_{n-2},
\end{equation*}
and
\begin{equation*}
v_{n+1} = \frac{(a+n) (b+n)-i p (c+2 n)}{(n+1) (c+n)} v_n + \frac{p (p+i (a+b+2 n-1))}{(n+1) (c+n)} v_{n-1} -\frac{p^2}{(n+1) (c+n)} v_{n-2}.
\end{equation*}
\end{thm}

\begin{thm}  \label{thm-cos-F-2}
Let $a,b,c,p,z \in \mathbb{C}$ and $-c \notin \mathbb{N} \cup \{0\}$. Suppose that
\begin{equation*}
\cos(pz) F(a,b;c;z) = \sum_{n=0}^\infty \frac{u_n + v_n}{2} z^n, \quad |z|<1.
\end{equation*}
Then $u_0 = 1$, $u_1 = \frac{ab}{c}+ip$, $u_2=\frac{i a b p}{c}+\frac{a (a+1) b (b+1)}{2 c (c+1)}-\frac{p^2}{2}$, $v_0 = 1$, $v_1 = \frac{ab}{c}-ip$, $u_2=-\frac{i a b p}{c}+\frac{a (a+1) b (b+1)}{2 c (c+1)}-\frac{p^2}{2}$
\begin{equation*}
u_{n+1} = \frac{(a+n) (b+n)+i p (c+2 n)}{(n+1) (c+n)} u_n + \frac{p (p-i (a+b+2 n-1))}{(n+1) (c+n)} u_{n-1} -\frac{p^2}{(n+1) (c+n)} u_{n-2},
\end{equation*}
and
\begin{equation*}
v_{n+1} = \frac{(a+n) (b+n)-i p (c+2 n)}{(n+1) (c+n)} v_n + \frac{p (p+i (a+b+2 n-1))}{(n+1) (c+n)} v_{n-1} -\frac{p^2}{(n+1) (c+n)} v_{n-2}.
\end{equation*}
\end{thm}

Next, we consider the product of $(1-\theta z)^p$ and the Gauss hypergeometric function. The recurrence relations for the coefficients of this function in the real domain were studied in \cite{Yang-PAMS-2025}, and were subsequently extended to the complex domain by Xu et al.\ \cite{Xu-Arxiv-2026}.

\begin{thm} \label{thm-theta-p-F}
Let $a,b,c,\theta,p,z \in \mathbb{C}$ and $-c \notin \mathbb{N} \cup \{0\}$. Suppose that
\begin{equation*}
(1-\theta z)^p F(a,b;c;z) = \sum_{n=0}^\infty u_n z^n, \quad |z| < \frac{1}{|\theta|}.
\end{equation*}
Then $u_0 = 1$, $u_1 = \frac{ab}{c}-\theta  p$, $u_2 = -\frac{a b \theta  p}{c}+\frac{a (a+1) b (b+1)}{2 c (c+1)}+\frac{1}{2} \theta ^2 (p-1) p$ and
\begin{equation*}
u_{n+1} = \alpha_0(n) u_n + \alpha_1(n) u_{n-1} + \alpha_2(n) u_{n-2},
\end{equation*}
where
\begin{equation*}
\begin{aligned}
\alpha_0(n) &= \frac{(a+n) (b+n)+2 \theta  n (c+n-1)-\theta  p (c+2 n)}{(n+1) (c+n)},\\
\alpha_1(n) &= \frac{\theta  (a (p-2 (b+n-1))+b (-2 n+p+2)+(c-2) \theta -(n-p) (\theta  (c+n-p-3)+2 n)+4 n-p-2)}{(n+1) (c+n)},\\
\alpha_2(n) &= \frac{\theta ^2 (a+n-p-2) (b+n-p-2)}{(n+1) (c+n)}.
\end{aligned}
\end{equation*}
\end{thm}

We now consider the product of $e^{-p \arctan z}$ and the Gauss hypergeometric function.

\begin{thm} \label{thm-e-arctan-F}
Let $a,b,c,p,z \in \mathbb{C}$ and $-c \notin \mathbb{N} \cup \{0\}$. Suppose that
\begin{equation*}
e^{-p \arctan z} F(a,b;c;z) = \sum_{n=0}^\infty u_n z^n, \quad |z|<1.
\end{equation*}
Then $u_0 = 1$, $u_1 = \frac{ab}{c}- p$, $u_2 = -\frac{a b p}{c}+\frac{a (a+1) b (b+1)}{2 c (c+1)}+\frac{p^2}{2}$, $u_3 = \frac{a b p^2}{2 c}-\frac{a (a+1) b (b+1) p}{2 c (c+1)}+\frac{a (a+1) (a+2) b (b+1) (b+2)}{6 c (c+1) (c+2)}+\frac{1}{3} \left(p-\frac{p^3}{2}\right)$,
\begin{equation*}
u_4 = \frac{\left(
\begin{aligned}
& \frac{6 a (a+1) b (b+1) p^2}{c (c+1)}-\frac{4 a b \left(p^2-2\right) p}{c}-\frac{4 a (a+1) (a+2) b (b+1) (b+2) p}{c (c+1) (c+2)} \\
& \quad +\frac{a (a+1) (a+2) (a+3) b (b+1) (b+2) (b+3)}{c (c+1) (c+2) (c+3)}+p^4-8 p^2
\end{aligned}
\right)}{24}
\end{equation*}
and
\begin{equation*}
u_{n+1} = \beta_0(n) u_n + \beta_1(n) u_{n-1} + \beta_2(n) u_{n-2} + \beta_3(n) u_{n-3} + \beta_4(n) u_{n-4},
\end{equation*}
where
\begin{equation*}
\begin{aligned}
\beta_0(n) &= \frac{(a+n) (b+n)-p (c+2 n)}{(n+1) (c+n)},\\
\beta_1(n) &= \frac{p (a+b+2 n-1)-2 (n-1) (c+n-2)-p^2}{(n+1) (c+n)},\\
\beta_2(n) &= \frac{2 (a+n-2) (b+n-2)-p (c+2 n-6)+p^2}{(n+1) (c+n)}, \\
\beta_3(n) &= \frac{p (a+b+2 n-7)-(n-3) (c+n-4)}{(n+1) (c+n)}, \\
\beta_4(n) &= \frac{(a+n-4) (b+n-4)}{(n+1) (c+n)}.
\end{aligned}
\end{equation*}
\end{thm}

In Theorems~\ref{thm-sin-F-2} and~\ref{thm-cos-F-2}, we have investigated the recurrence relations satisfied by the coefficients of the products of trigonometric functions and the Gauss hypergeometric function, which are characterized as linear combinations of two recurrence relations. The following two theorems now provide a single recurrence relation satisfied by these coefficients.

\begin{thm} \label{thm-sin-F}
Let $a,b,c,p,z \in \mathbb{C}$ and $-c \notin \mathbb{N} \cup \{0\}$. Suppose that
\begin{equation*}
\sin(pz) F(a,b;c;z) = \sum_{n=0}^\infty u_n z^n, \quad |z|<1.
\end{equation*}
Then
\begin{equation*}
u_0 = 0, \quad u_1 = p, \quad u_2 = \frac{a b p}{c}, \quad u_3 = \frac{a (a+1) b (b+1) p}{2 c (c+1)}-\frac{p^3}{6},
\end{equation*}
\begin{equation*}
u_4 = \frac{a (a+1) (a+2) b (b+1) (b+2) p}{6 c (c+1) (c+2)}-\frac{a b p^3}{6 c},
\end{equation*}
\begin{equation*}
u_5 = -\frac{a (a+1) b (b+1) p^3}{12 c (c+1)}+\frac{a (a+1) (a+2) (a+3) b (b+1) (b+2) (b+3) p}{24 c (c+1) (c+2) (c+3)}+\frac{p^5}{120},
\end{equation*}
\begin{equation*}
u_6 = \frac{a b p^5}{120 c}-\frac{a (a+1) (a+2) b (b+1) (b+2) p^3}{36 c (c+1) (c+2)}+\frac{a (a+1) (a+2) (a+3) (a+4) b (b+1) (b+2) (b+3) (b+4) p}{120 c (c+1) (c+2) (c+3) (c+4)},
\end{equation*}
\begin{equation*}
\begin{aligned}
u_7 & = \frac{a (a+1) b (b+1) p^5}{240 c (c+1)}-\frac{a (a+1) (a+2) (a+3) b (b+1) (b+2) (b+3) p^3}{144 c (c+1) (c+2) (c+3)} \\
& \quad +\frac{a (a+1) (a+2) (a+3) (a+4) (a+5) b (b+1) (b+2) (b+3) (b+4) (b+5) p}{720 c (c+1) (c+2) (c+3) (c+4) (c+5)}-\frac{p^7}{5040},
\end{aligned}
\end{equation*}
\begin{equation*}
\begin{aligned}
u_8 &= -\frac{a b p^7}{5040 c}+\frac{a (a+1) (a+2) b (b+1) (b+2) p^5}{720 c (c+1) (c+2)} \\
& \quad -\frac{a (a+1) (a+2) (a+3) (a+4) b (b+1) (b+2) (b+3) (b+4) p^3}{720 c (c+1) (c+2) (c+3) (c+4)} \\
& \quad +\frac{a (a+1) (a+2) (a+3) (a+4) (a+5) (a+6) b (b+1) (b+2) (b+3) (b+4) (b+5) (b+6) p}{5040 c (c+1) (c+2) (c+3) (c+4) (c+5) (c+6)},
\end{aligned}
\end{equation*}
\begin{equation*}
\begin{aligned}
u_9 &= -\frac{a (a+1) b (b+1) p^7}{10080 c (c+1)}+ \frac{a (a+1) (a+2) (a+3) b (b+1) (b+2) (b+3) p^5}{2880 c (c+1) (c+2) (c+3)} \\
& \quad -\frac{a (a+1) (a+2) (a+3) (a+4) (a+5) b (b+1) (b+2) (b+3) (b+4) (b+5) p^3}{4320 c (c+1) (c+2) (c+3) (c+4) (c+5)} \\
& \quad +\frac{a (a+1) (a+2) (a+3) (a+4) (a+5) (a+6) (a+7) b (b+1) (b+2) (b+3) (b+4) (b+5) (b+6) (b+7) p}{40320 c (c+1) (c+2) (c+3) (c+4) (c+5) (c+6) (c+7)} \\
& \quad +\frac{p^9}{362880},
\end{aligned}
\end{equation*}
and
\begin{equation*}
u_{n+1} = \sum_{i=0}^9 \gamma_i(n) u_{n-i}, \quad n \geq 9,
\end{equation*}
where
\begin{equation*}
\gamma_0(n) = \frac{2 (n-1)^2 (a (c-2 b)+c (b+c-3))+4 (c-2) (n-1) (c (a+b)-a b)+2 a b (c-2) (c+1)}{(c-2) c (n+1) (c+n)},
\end{equation*}
\begin{equation*}
\gamma_1(n) = \frac{
\left(
\begin{aligned}
& -(n-4) (n-3) (n-2) (n-1) \left(a^2+4 c (a+b)-10 a b+b^2+c^2-6 c+4 p^2-1\right) \\
& \quad +2 (n-3) (n-2) (n-1) \left(
\begin{aligned}
& a^2 (4 b-3 c)+a \left(2 (b-1) c+4 b (b+3)-3 c^2\right) \\
& \quad -c \left(b (3 b+3 c+2)+c+4 p^2-13\right)
\end{aligned}\right) \\
& \quad +(n-2) (n-1) \left(
\begin{aligned}
& a^2 \left(8 b^2+b (4 c+6)-6 c^2+c\right) \\
& \quad +a \left(-(10 b+7) c^2+4 (b (b+3)+3) c+6 b (b+1)\right) \\
& \quad +c \left(b^2 (1-6 c)+b (12-7 c)-6 (c-2) p^2+c+11\right)
\end{aligned}\right) \\
& \quad +2 (n-1) \left(a^2 b (b (4 c-2)-3 (c-1) c)-a b c (3 b (c-1)+c-5)+c \left(-c^2+c+2\right) p^2\right) \\
& \quad -(c-2) \left(a^2 b (b (c+2)-c)-a b (b+1) c+c^2 (c+1) p^2\right)
\end{aligned}
\right)
}{(c-2) c n (n+1) (c+n-1) (c+n)},
\end{equation*}
\begin{equation*}
\gamma_2(n) = \frac{2\left(
\begin{aligned}
& (n-5) (n-4) (n-3) (n-2) \left(a^2+a (c-4 b)+(b-1) (b+c+1)+8 p^2\right) \\
& \quad +(n-4) (n-3) (n-2) \left(
\begin{aligned}
& a^3+a^2 (-5 b+3 c+2)+a \left(-5 b^2+2 b (c-7)+3 c+4 p^2-1\right) \\
& \quad +(b+3 c+1) \left(b^2+b+4 p^2-2\right)
\end{aligned}\right) \\
& \quad -(n-3) (n-2) \left(
\begin{aligned}
& a^3 (b-2 c)+a^2 (b (10 b+9)-4 (b+1) c) \\
& \quad +a \left(b (b+1) (b-4 c+8)-6 c p^2+2 c\right)\\
& \quad +2 c \left(-b^3-2 b^2-3 p^2 (b+c-3)+b+2\right)
\end{aligned}\right) \\
& \quad -(n-2) \left(
\begin{aligned}
& a^3 b (4 b-3 c+2)+a^2 b (2 b+1) (2 b-c)+a p^2 \left(10 b-3 c^2+c\right) \\
& \quad -a (b-1) b (b (3 c-2)+4 c-2)+c p^2 \left(-3 b c+b-c^2+9\right)
\end{aligned}\right) \\
& \quad +(c+1) p^2 \left(c^2 (2 a+2 b+1)-3 (a+1) (b+1) c+4 a b\right) \\
& \quad -(a-1) a (b-1) b (-c (a+b+1)+2 a b+a+b+1)
\end{aligned}
\right)}{(c-2) c n (n+1) (c+n-1) (c+n)},
\end{equation*}
\begin{equation*}
\gamma_3(n) = -\frac{\left(
\begin{aligned}
& (n-6) (n-5) (n-4) (n-3) \left((a-b)^2+24 p^2-1\right) \\
& \quad +2 (n-5) (n-4) (n-3) \left(
\begin{aligned}
& a^3-a^2 (b-2)-a \left(b (b+4)-12 p^2+1\right) \\
& \quad +b^3+2 b^2+12 p^2 (b+c+1)-b-2
\end{aligned}\right) \\
& \quad +(n-4) (n-3) \left(
\begin{aligned}
& a^4+a^3 (2 b+3)+a^2 \left(-3 b (2 b+1)+6 p^2+1\right) \\
& \quad +a \left(12 p^2 (b+2 c)+b (b (2 b-3)-8)-3\right) \\
& \quad +6 p^2 \left(b^2+4 b c+(c-6) c-1\right)+(b+1)^2 \left(b^2+b-2\right)+8 p^4
\end{aligned}\right) \\
& \quad  +2 (n-3) \left(
\begin{aligned}
& a^4 b-a^3 b^2-a^2 \left(b^3+b-3 c p^2\right)+a b^2 \left(b^2-1\right) \\
& \quad +a p^2 (b (6 c-40)+c (3 c+2))+c p^2 \left(b (3 b+3 c+2)+c+4 p^2-13\right)
\end{aligned}\right) \\
& \quad +a^4 (b-1) b+a^3 b \left(-2 b^2+b+1\right)\\
& \quad+a^2 \left(b^4+b^3+p^2 \left(12 b^2-2 b (6 c+5)+c (6 c-1)\right)-3 b^2+b\right) \\
& \quad +a p^2 \left((6 b+7) c^2-4 (b (3 b+2)+3) c+2 b (11-5 b)\right)-a (b-1)^2 b (b+1)\\
& \quad+c p^2 \left(b^2 (6 c-1)+b (7 c-12)+5 (c-2) p^2-c-11\right)
\end{aligned}
\right)}{(c-2) c n (n+1) (c+n-1) (c+n)},
\end{equation*}
\begin{equation*}
\gamma_4(n) = \frac{\left(
\begin{aligned}
& 8 (n-7) (n-6) (n-5) (n-4) +4 (n-6) (n-5) (n-4) (3 (a+b+1)+c) \\
& \quad +2 (n-5) (n-4) \left(3 (a+b-1) (a+b+c+1)+8 p^2\right) \\
& \quad +(n-4) \left(
\begin{aligned}
& a^3+a^2 (3 b+3 c+2)+a \left(b (3 b+6 c-46)+3 c+4 p^2-1\right)\\
& \quad +(b+3 c+1) \left(b^2+b+4 p^2-2\right)
\end{aligned}\right) \\
& \quad 2 p^2 \left(
\begin{aligned}
& a^3 (2 c-3 b)+a^2 (3 b (2 b-5)+4 c)-a \left(p^2 (6 b-5 c)-4 b c+3 b (b (b+5)-4)+2 c\right) \\
& \quad +5 c p^2 (b+c-3)+2 (b-1) (b+1) (b+2) c
\end{aligned}\right)
\end{aligned}
\right)}{(c-2) c n (n+1) (c+n-1) (c+n)},
\end{equation*}
\begin{equation*}
\gamma_5(n) = \frac{\left(
\begin{aligned}
& -4 (n-8) (n-7) (n-6) (n-5) -8 (n-7) (n-6) (n-5) (a+b+1) \\
& \quad -6 (n-6) (n-5) \left((a+b)^2+8 p^2-1\right) \\
& \quad -2 (n-5) \left(a^3+3 a^2 b+2 a^2+3 a b^2+12 p^2 (a+b+c+1)-16 a b-a+b^3+2 b^2-b-2\right) \\
& \quad p^2 \left(-a^4+a^3 (2 b-3)-a^2 \left(b (6 b-11)+5 p^2+1\right)+a \left(2 p^2 (13 b-10 c)+b (b (2 b+11)-12)+3\right) \right) \\
& \quad + p^2  \left( -5 p^2 \left(b^2+4 b c+(c-6) c-1\right)-(b+1)^2 \left(b^2+b-2\right)-4 p^4\right)
\end{aligned}
\right)}{(c-2) c n (n+1) (c+n-1) (c+n)},
\end{equation*}
\begin{equation*}
\gamma_6(n) = \frac{2 p^4 \left(5 a^2+a (-8 b+5 c+12 n-72)+5 b^2+5 b c+12 b (n-6)+4 n (c+4 n)-29 c-196 n+8 p^2+595\right)}{(c-2) c n (n+1) (c+n-1) (c+n)},
\end{equation*}
\begin{equation*}
\gamma_7(n) = -\frac{p^4 \left(5 a^2-2 a (b-4 n+28)+5 b^2+8 b (n-7)+8 (n-14) n+24 p^2+387\right)}{(c-2) c n (n+1) (c+n-1) (c+n)},
\end{equation*}
\begin{equation*}
\gamma_8(n) = \frac{16 p^6}{(c-2) c n (n+1) (c+n-1) (c+n)},
\end{equation*}
\begin{equation*}
\gamma_9(n) = -\frac{4 p^6}{(c-2) c n (n+1) (c+n-1) (c+n)}.
\end{equation*}
\end{thm}

\begin{thm} \label{thm-cos-F}
Let $a,b,c,p,z \in \mathbb{C}$ and $-c \notin \mathbb{N} \cup \{0\}$. Suppose that
\begin{equation*}
\cos(pz) F(a,b;c;z) = \sum_{n=0}^\infty u_n z^n, \quad |z|<1.
\end{equation*}
Then
\begin{equation*}
u_0 = 1, \quad u_1 = \frac{a b}{c}, \quad u_2 = \frac{a (a+1) b (b+1)}{2 c (c+1)}-\frac{p^2}{2}, \quad u_3 = \frac{a (a+1) (a+2) b (b+1) (b+2)}{6 c (c+1) (c+2)}-\frac{a b p^2}{2 c},
\end{equation*}
\begin{equation*}
u_4 = -\frac{a (a+1) b (b+1) p^2}{4 c (c+1)}+\frac{a (a+1) (a+2) (a+3) b (b+1) (b+2) (b+3)}{24 c (c+1) (c+2) (c+3)}+\frac{p^4}{24},
\end{equation*}
\begin{equation*}
u_5 = \frac{a b p^4}{24 c}-\frac{a (a+1) (a+2) b (b+1) (b+2) p^2}{12 c (c+1) (c+2)}+\frac{a (a+1) (a+2) (a+3) (a+4) b (b+1) (b+2) (b+3) (b+4)}{120 c (c+1) (c+2) (c+3) (c+4)},
\end{equation*}
\begin{equation*}
\begin{aligned}
u_6 & = \frac{a (a+1) b (b+1) p^4}{48 c (c+1)}-\frac{a (a+1) (a+2) (a+3) b (b+1) (b+2) (b+3) p^2}{48 c (c+1) (c+2) (c+3)} \\
& \quad +\frac{a (a+1) (a+2) (a+3) (a+4) (a+5) b (b+1) (b+2) (b+3) (b+4) (b+5)}{720 c (c+1) (c+2) (c+3) (c+4) (c+5)}-\frac{p^6}{720},
\end{aligned}
\end{equation*}
\begin{equation*}
\begin{aligned}
u_7 &= -\frac{a b p^6}{720 c}+\frac{a (a+1) (a+2) b (b+1) (b+2) p^4}{144 c (c+1) (c+2)}-\frac{a (a+1) (a+2) (a+3) (a+4) b (b+1) (b+2) (b+3) (b+4) p^2}{240 c (c+1) (c+2) (c+3) (c+4)} \\
& \quad +\frac{a (a+1) (a+2) (a+3) (a+4) (a+5) (a+6) b (b+1) (b+2) (b+3) (b+4) (b+5) (b+6)}{5040 c (c+1) (c+2) (c+3) (c+4) (c+5) (c+6)},
\end{aligned}
\end{equation*}
\begin{equation*}
\begin{aligned}
u_8 &= -\frac{a (a+1) b (b+1) p^6}{1440 c (c+1)}+\frac{a (a+1) (a+2) (a+3) b (b+1) (b+2) (b+3) p^4}{576 c (c+1) (c+2) (c+3)} \\
& \quad -\frac{a (a+1) (a+2) (a+3) (a+4) (a+5) b (b+1) (b+2) (b+3) (b+4) (b+5) p^2}{1440 c (c+1) (c+2) (c+3) (c+4) (c+5)} \\
& \quad +\frac{a (a+1) (a+2) (a+3) (a+4) (a+5) (a+6) (a+7) b (b+1) (b+2) (b+3) (b+4) (b+5) (b+6) (b+7)}{40320 c (c+1) (c+2) (c+3) (c+4) (c+5) (c+6) (c+7)} \\
& \quad +\frac{p^8}{40320},
\end{aligned}
\end{equation*}
\begin{equation*}
\begin{aligned}
u_9
& = \frac{a b p^8}{40320 c}-\frac{a (a+1) (a+2) b (b+1) (b+2) p^6}{4320 c (c+1) (c+2)} \\
& \quad +\frac{a (a+1) (a+2) (a+3) (a+4) b (b+1) (b+2) (b+3) (b+4) p^4}{2880 c (c+1) (c+2) (c+3) (c+4)} \\
& \quad -\frac{a (a+1) (a+2) (a+3) (a+4) (a+5) (a+6) b (b+1) (b+2) (b+3) (b+4) (b+5) (b+6) p^2}{10080 c (c+1) (c+2) (c+3) (c+4) (c+5) (c+6)} \\
& \quad +\frac{ (a)_9 (b)_9 }{362880 (c)_9},
\end{aligned}
\end{equation*}
and
\begin{equation*}
u_{n+1} = \sum_{i=0}^9 \gamma_i(n) u_{n-i}, \quad n \geq 9,
\end{equation*}
where $\gamma_i(n)$ ($i=0,1,\cdots,9$) are defined in Theorem \ref{thm-sin-F}.
\end{thm}

Similar, in Theorems~\ref{thm-sinh-F-2} and~\ref{thm-cosh-F-2}, we have investigated the recurrence relations satisfied by the coefficients of the products of hyperbolic functions and the Gauss hypergeometric function, which are characterized as linear combinations of two recurrence relations. The following two theorems now provide a single recurrence relation satisfied by these coefficients.

\begin{thm} \label{thm-sinh-F}
Let $a,b,c,p,z \in \mathbb{C}$ and $-c \notin \mathbb{N} \cup \{0\}$. Suppose that
\begin{equation*}
\sinh(pz) F(a,b;c;z) = \sum_{n=0}^\infty u_n z^n, \quad |z|<1.
\end{equation*}
Then
\begin{equation*}
u_0 = 0, \quad u_1 = p, \quad u_2 = \frac{a b p}{c}, \quad u_3 = \frac{a (a+1) b (b+1) p}{2 c (c+1)}+\frac{p^3}{6}, \quad u_4 = \frac{a b p^3}{6 c}+\frac{a (a+1) (a+2) b (b+1) (b+2) p}{6 c (c+1) (c+2)},
\end{equation*}
\begin{equation*}
u_5 = \frac{a (a+1) b (b+1) p^3}{12 c (c+1)}+\frac{a (a+1) (a+2) (a+3) b (b+1) (b+2) (b+3) p}{24 c (c+1) (c+2) (c+3)}+\frac{p^5}{120},
\end{equation*}
\begin{equation*}
u_6 = \frac{a b p^5}{120 c}+\frac{a (a+1) (a+2) b (b+1) (b+2) p^3}{36 c (c+1) (c+2)}+\frac{a (a+1) (a+2) (a+3) (a+4) b (b+1) (b+2) (b+3) (b+4) p}{120 c (c+1) (c+2) (c+3) (c+4)},
\end{equation*}
\begin{equation*}
\begin{aligned}
u_7 &= \frac{a (a+1) b (b+1) p^5}{240 c (c+1)}+\frac{a (a+1) (a+2) (a+3) b (b+1) (b+2) (b+3) p^3}{144 c (c+1) (c+2) (c+3)}\\
& \quad +\frac{a (a+1) (a+2) (a+3) (a+4) (a+5) b (b+1) (b+2) (b+3) (b+4) (b+5) p}{720 c (c+1) (c+2) (c+3) (c+4) (c+5)}+\frac{p^7}{5040},
\end{aligned}
\end{equation*}
\begin{equation*}
\begin{aligned}
u_8 &= \frac{a b p^7}{5040 c}+\frac{a (a+1) (a+2) b (b+1) (b+2) p^5}{720 c (c+1) (c+2)}+\frac{a (a+1) (a+2) (a+3) (a+4) b (b+1) (b+2) (b+3) (b+4) p^3}{720 c (c+1) (c+2) (c+3) (c+4)} \\
& \quad +\frac{a (a+1) (a+2) (a+3) (a+4) (a+5) (a+6) b (b+1) (b+2) (b+3) (b+4) (b+5) (b+6) p}{5040 c (c+1) (c+2) (c+3) (c+4) (c+5) (c+6)},
\end{aligned}
\end{equation*}
\begin{equation*}
\begin{aligned}
u_9 &= \frac{a (a+1) b (b+1) p^7}{10080 c (c+1)}+\frac{a (a+1) (a+2) (a+3) b (b+1) (b+2) (b+3) p^5}{2880 c (c+1) (c+2) (c+3)}\\
& \quad +\frac{a (a+1) (a+2) (a+3) (a+4) (a+5) b (b+1) (b+2) (b+3) (b+4) (b+5) p^3}{4320 c (c+1) (c+2) (c+3) (c+4) (c+5)} \\
& \quad +\frac{a (a+1) (a+2) (a+3) (a+4) (a+5) (a+6) (a+7) b (b+1) (b+2) (b+3) (b+4) (b+5) (b+6) (b+7) p}{40320 c (c+1) (c+2) (c+3) (c+4) (c+5) (c+6) (c+7)}\\
& \quad+\frac{p^9}{362880},
\end{aligned}
\end{equation*}
and
\begin{equation*}
u_{n+1} = \sum_{i=0}^9 \delta_i(n) u_{n-i}, \quad n \geq 9,
\end{equation*}
where
\begin{equation*}
\delta_0(n) = \frac{2 (n-1)^2 (a (c-2 b)+c (b+c-3))+4 (c-2) (n-1) (c (a+b)-a b)+2 a b (c-2) (c+1)}{(c-2) c (n+1) (c+n)},
\end{equation*}
\begin{equation*}
\delta_1(n) = \frac{
\left(
\begin{aligned}
& -(n-4) (n-3) (n-2) (n-1) \left(a^2+4 c (a+b)-10 a b+b^2+c^2-6 c-4 p^2-1\right) \\
& \quad +2 (n-3) (n-2) (n-1) \left(
\begin{aligned}
& a^2 (4 b-3 c)+a \left(2 (b-1) c+4 b (b+3)-3 c^2\right) \\
& \quad -c \left(b (3 b+3 c+2)+c-4 p^2-13\right)
\end{aligned}\right) \\
& \quad +(n-2) (n-1) \left(
\begin{aligned}
& a^2 \left(8 b^2+b (4 c+6)-6 c^2+c\right) \\
& \quad +a \left(-(10 b+7) c^2+4 (b (b+3)+3) c+6 b (b+1)\right) \\
& \quad +c \left(b^2 (1-6 c)+b (12-7 c)+6 (c-2) p^2+c+11\right)
\end{aligned}\right) \\
& \quad +(c-2) \left(a^2 b (c-b (c+2))+a b (b+1) c+c^2 (c+1) p^2\right)\\
& \quad +2 (n-1) \left(a^2 b (b (4 c-2)-3 (c-1) c)-a b c (3 b (c-1)+c-5)+(c-2) c (c+1) p^2\right)
\end{aligned}
\right)
}{(c-2) c n (n+1) (c+n-1) (c+n)},
\end{equation*}
\begin{equation*}
\delta_2(n) = \frac{2\left(
\begin{aligned}
& (n-5) (n-4) (n-3) (n-2) \left(a^2+a (c-4 b)+(b-1) (b+c+1)-8 p^2\right) \\
& \quad + (n-4) (n-3) (n-2) \left(
\begin{aligned}
& a^3+a^2 (-5 b+3 c+2)-a \left(b (5 b-2 c+14)-3 c+4 p^2+1\right) \\
& \quad +(b+3 c+1) \left(b^2+b-4 p^2-2\right)
\end{aligned}\right) \\
& \quad - (n-3) (n-2) \left(
\begin{aligned}
&a^3 (b-2 c)+a^2 (b (10 b+9)-4 (b+1) c)\\
& \quad +a \left(b (b+1) (b-4 c+8)+6 c p^2+2 c\right)\\
& \quad +2 c \left(-b^3-2 b^2+3 p^2 (b+c-3)+b+2\right)
\end{aligned}\right)\\
& \quad - (n-2) \left(
\begin{aligned}
& a^3 b (4 b-3 c+2)+a^2 b (2 b+1) (2 b-c) \\
& \quad -a \left(p^2 \left(10 b-3 c^2+c\right)+(b-1) b (b (3 c-2)+4 c-2)\right) \\
& \quad +c p^2 \left(b (3 c-1)+c^2-9\right)
\end{aligned}\right)\\
& \quad - (c+1) p^2 \left(c^2 (2 a+2 b+1)-3 (a+1) (b+1) c+4 a b\right) \\
& \quad -(a-1) a (b-1) b (-c (a+b+1)+2 a b+a+b+1)
\end{aligned}
\right)}{(c-2) c n (n+1) (c+n-1) (c+n)},
\end{equation*}
\begin{equation*}
\delta_3(n) = -\frac{\left(
\begin{aligned}
& (n-6) (n-5) (n-4) (n-3) \left((a-b)^2-24 p^2-1\right) \\
& \quad +2 (n-5) (n-4) (n-3) \left(
\begin{aligned}
& a^3-a^2 (b-2)-a \left(b (b+4)+12 p^2+1\right) \\
& \quad +b^3+2 b^2-12 p^2 (b+c+1)-b-2
\end{aligned}\right) \\
& \quad +(n-4) (n-3) \left(
\begin{aligned}
& a^4+a^3 (2 b+3)-a^2 \left(6 b^2+3 b+6 p^2-1\right) \\
& \quad +a \left(-12 p^2 (b+2 c)+b (b (2 b-3)-8)-3\right) \\
& \quad -6 p^2 \left(b^2+4 b c+(c-6) c-1\right)+(b+1)^2 \left(b^2+b-2\right)+8 p^4
\end{aligned}
\right)\\
& \quad +2 (n-3) \left(
\begin{aligned}
& -p^2 \left(c^2 (3 a+3 b+1)+c (a+b) (3 a+3 b+2)-40 a b-13 c\right) \\
& \quad +a b (a-b-1) (a-b+1) (a+b)+4 c p^4
\end{aligned}\right) \\
& \quad + a^4 (b-1) b+a^3 b \left(-2 b^2+b+1\right)+a^2 \left(b^4+b^3+p^2 \left(-12 b^2+2 b (6 c+5)-6 c^2+c\right)-3 b^2+b\right) \\
& \quad +a p^2 \left(-(6 b+7) c^2+4 (b (3 b+2)+3) c+2 b (5 b-11)\right) \\
& \quad -a (b-1)^2 b (b+1)+c p^2 \left(b^2 (1-6 c)+b (12-7 c)+5 (c-2) p^2+c+11\right)
\end{aligned}
\right)}{(c-2) c n (n+1) (c+n-1) (c+n)},
\end{equation*}
\begin{equation*}
\delta_4(n) = \frac{2 p^2\left(
\begin{aligned}
& -8 (n-7) (n-6) (n-5) (n-4) -4 (n-6) (n-5) (n-4) (3 (a+b+1)+c) \\
& \quad +2 (n-5) (n-4) \left(8 p^2-3 (a+b-1) (a+b+c+1)\right) \\
& \quad +(n-4) \left(
\begin{aligned}
& -a^3-a^2 (3 b+3 c+2)+a \left(b (-3 b-6 c+46)-3 c+4 p^2+1\right) \\
& \quad -(b+3 c+1) \left(b^2+b-4 p^2-2\right)
\end{aligned}\right) \\
& \quad a^3 (3 b-2 c)+a^2 (3 (5-2 b) b-4 c)+a \left(p^2 (5 c-6 b)-4 b c+3 b (b (b+5)-4)+2 c\right) \\
& \quad +5 c p^2 (b+c-3)-2 (b-1) (b+1) (b+2) c
\end{aligned}
\right)}{(c-2) c n (n+1) (c+n-1) (c+n)},
\end{equation*}
\begin{equation*}
\delta_5(n) = \frac{p^2 \left(
\begin{aligned}
& 4 (n-8) (n-7) (n-6) (n-5) +8 (n-7) (n-6) (n-5) (a+b+1) \\
& \quad +6 (n-6) (n-5) \left((a+b)^2-8 p^2-1\right) \\
& \quad -2 (n-5) \left(-a^3-3 a^2 b-2 a^2-3 a b^2+12 p^2 (a+b+c+1)+16 a b+a-b^3-2 b^2+b+2\right) \\
& \quad + a^4+a^3 (3-2 b)+a^2 \left(b (6 b-11)-5 p^2+1\right)-a \left(2 b^3+11 b^2-2 b \left(13 p^2+6\right)+20 c p^2+3\right) \\
& -5 p^2 \left(b^2+4 b c+(c-6) c-1\right)+(b+1)^2 \left(b^2+b-2\right)+4 p^4
\end{aligned}
\right)}{(c-2) c n (n+1) (c+n-1) (c+n)},
\end{equation*}
\begin{equation*}
\delta_6(n) = \frac{2 p^4 \left(5 a^2+a (-8 b+5 c+12 n-72)+5 b^2+5 b c+12 b (n-6)+4 n (c+4 n)-29 c-196 n-8 p^2+595\right)}{(c-2) c n (n+1) (c+n-1) (c+n)},
\end{equation*}
\begin{equation*}
\delta_7(n) = \frac{p^4 \left(-5 a^2+2 a (b-4 n+28)-5 b^2-8 b (n-7)-8 (n-14) n+24 p^2-387\right)}{(c-2) c n (n+1) (c+n-1) (c+n)},
\end{equation*}
\begin{equation*}
\delta_8(n) = -\frac{16 p^6}{(c-2) c n (n+1) (c+n-1) (c+n)},
\end{equation*}
\begin{equation*}
\delta_9(n) = \frac{4 p^6}{(c-2) c n (n+1) (c+n-1) (c+n)}.
\end{equation*}
\end{thm}

\begin{thm} \label{thm-cosh-F}
Let $a,b,c,p,z \in \mathbb{C}$ and $-c \notin \mathbb{N} \cup \{0\}$. Suppose that
\begin{equation*}
\cosh(pz) F(a,b;c;z) = \sum_{n=0}^\infty u_n z^n, \quad |z|<1.
\end{equation*}
Then
\begin{equation*}
u_0 = 1, \quad u_1 = \frac{a b}{c}, \quad u_2 = \frac{a (a+1) b (b+1)}{2 c (c+1)}+\frac{p^2}{2}, \quad u_3 = \frac{a b p^2}{2 c}+\frac{a (a+1) (a+2) b (b+1) (b+2)}{6 c (c+1) (c+2)},
\end{equation*}
\begin{equation*}
u_4 = \frac{a (a+1) b (b+1) p^2}{4 c (c+1)}+\frac{a (a+1) (a+2) (a+3) b (b+1) (b+2) (b+3)}{24 c (c+1) (c+2) (c+3)}+\frac{p^4}{24},
\end{equation*}
\begin{equation*}
u_5 = \frac{a b p^4}{24 c}+\frac{a (a+1) (a+2) b (b+1) (b+2) p^2}{12 c (c+1) (c+2)}+\frac{a (a+1) (a+2) (a+3) (a+4) b (b+1) (b+2) (b+3) (b+4)}{120 c (c+1) (c+2) (c+3) (c+4)},
\end{equation*}
\begin{equation*}
\begin{aligned}
u_6 & = \frac{a (a+1) b (b+1) p^4}{48 c (c+1)}+\frac{a (a+1) (a+2) (a+3) b (b+1) (b+2) (b+3) p^2}{48 c (c+1) (c+2) (c+3)} \\
& \quad +\frac{a (a+1) (a+2) (a+3) (a+4) (a+5) b (b+1) (b+2) (b+3) (b+4) (b+5)}{720 c (c+1) (c+2) (c+3) (c+4) (c+5)}+\frac{p^6}{720},
\end{aligned}
\end{equation*}
\begin{equation*}
\begin{aligned}
u_7 &= \frac{a b p^6}{720 c}+\frac{a (a+1) (a+2) b (b+1) (b+2) p^4}{144 c (c+1) (c+2)}+\frac{a (a+1) (a+2) (a+3) (a+4) b (b+1) (b+2) (b+3) (b+4) p^2}{240 c (c+1) (c+2) (c+3) (c+4)}\\
& \quad +\frac{a (a+1) (a+2) (a+3) (a+4) (a+5) (a+6) b (b+1) (b+2) (b+3) (b+4) (b+5) (b+6)}{5040 c (c+1) (c+2) (c+3) (c+4) (c+5) (c+6)},
\end{aligned}
\end{equation*}
\begin{equation*}
\begin{aligned}
u_8 &= \frac{a (a+1) b (b+1) p^6}{1440 c (c+1)}+\frac{a (a+1) (a+2) (a+3) b (b+1) (b+2) (b+3) p^4}{576 c (c+1) (c+2) (c+3)}\\
& \quad +\frac{a (a+1) (a+2) (a+3) (a+4) (a+5) b (b+1) (b+2) (b+3) (b+4) (b+5) p^2}{1440 c (c+1) (c+2) (c+3) (c+4) (c+5)} \\
& \quad +\frac{a (a+1) (a+2) (a+3) (a+4) (a+5) (a+6) (a+7) b (b+1) (b+2) (b+3) (b+4) (b+5) (b+6) (b+7)}{40320 c (c+1) (c+2) (c+3) (c+4) (c+5) (c+6) (c+7)}\\
& \quad +\frac{p^8}{40320}
\end{aligned}
\end{equation*}
\begin{equation*}
\begin{aligned}
u_9 &= \frac{a b p^8}{40320 c}+\frac{a (a+1) (a+2) b (b+1) (b+2) p^6}{4320 c (c+1) (c+2)}\\
& \quad +\frac{a (a+1) (a+2) (a+3) (a+4) b (b+1) (b+2) (b+3) (b+4) p^4}{2880 c (c+1) (c+2) (c+3) (c+4)} \\
& \quad +\frac{a (a+1) (a+2) (a+3) (a+4) (a+5) (a+6) b (b+1) (b+2) (b+3) (b+4) (b+5) (b+6) p^2}{10080 c (c+1) (c+2) (c+3) (c+4) (c+5) (c+6)}\\
& \quad +\frac{(a)_9 (b)_9}{362880 (c)_9}
\end{aligned}
\end{equation*}
and
\begin{equation*}
u_{n+1} = \sum_{i=0}^9 \delta_i(n) u_{n-i}, \quad n \geq 9,
\end{equation*}
where $\delta_i(n)$ ($i=0,1,\cdots,9$) are defined in Theorem \ref{thm-sinh-F}.
\end{thm}

\section{Recurrence relations for the Maclaurin coefficients of products of elementary functions and elliptic integrals}

Taking $(a,b,c)=(1/2,1/2,1)$ and $(a,b,c)=(-1/2,1/2,1)$ in Theorem \ref{thm-exp-F}, we obtain the following corollaries.

\begin{cor} \label{cor-exp-K}
Let $p,z \in \mathbb{C}$. Suppose that
\begin{equation*}
e^{pz} K(\sqrt{z}) = \frac{\pi}{2} \sum_{n=0}^\infty u_n z^n, \quad |z|<1.
\end{equation*}
Then $u_0 = 1$, $u_1 = \frac{1}{4} (4 p+1)$, $u_2=\frac{1}{64} \left(32 p^2+16 p+9\right)$ and
\begin{equation*}
u_{n+1} = \frac{(2 n+1) (2 n+4 p+1)}{4 (n+1)^2} u_n -\frac{p (2 n+p)}{(n+1)^2} u_{n-1} + \frac{p^2}{(n+1)^2} u_{n-2}.
\end{equation*}
\end{cor}

\begin{cor} \label{cor-exp-K}
Let $p,z \in \mathbb{C}$. Suppose that
\begin{equation*}
e^{pz} E(\sqrt{z}) = \frac{\pi}{2} \sum_{n=0}^\infty u_n z^n, \quad |z|<1.
\end{equation*}
Then $u_0 = 1$, $u_1 = \frac{1}{4} (4 p-1)$, $u_2=\frac{1}{64} \left(32 p^2-16 p-3\right)$ and
\begin{equation*}
u_{n+1} = \frac{(2 n+1) (2 n+4 p-1)}{4 (n+1)^2} u_n -\frac{p (2 n+p-1)}{(n+1)^2} u_{n-1} + \frac{p^2}{(n+1)^2} u_{n-2}.
\end{equation*}
\end{cor}

Taking $(a,b,c)=(1/2,1/2,1)$ and $(a,b,c)=(-1/2,1/2,1)$ in Theorem \ref{thm-theta-p-F}, we obtain the following corollaries.

\begin{cor}
Let $p, \theta, z \in \mathbb{C}$. Suppose that
\begin{equation*}
(1-\theta z)^p K(\sqrt{z}) = \frac{\pi}{2}  \sum_{n=0}^\infty u_n z^n, \quad |z| < \frac{1}{|\theta|}.
\end{equation*}
Then $u_0 = 1$, $u_1 = (1-4 \theta  p) /4$, $u_2 = \frac{1}{64} \left(32 \theta ^2 p^2-32 \theta ^2 p-16 \theta  p+9\right)$ and
\begin{equation*}
u_{n+1} = \alpha_0(n) u_n + \alpha_1(n) u_{n-1} + \alpha_2(n) u_{n-2},
\end{equation*}
where
\begin{equation*}
\begin{aligned}
\alpha_0(n) &= \frac{8 \theta  n (n-p)+(2 n+1)^2-4 \theta  p}{4 (n+1)^2},\\
\alpha_1(n) &= -\frac{\theta  \left(2 (\theta +2) n^2-4 (\theta +1) n (p+1)+2 \theta  (p+1)^2+1\right)}{2 (n+1)^2},\\
\alpha_2(n) &= \frac{\theta ^2 (-2 n+2 p+3)^2}{4 (n+1)^2}.
\end{aligned}
\end{equation*}
\end{cor}

\begin{cor}
Let $p, \theta, z \in \mathbb{C}$. Suppose that
\begin{equation*}
(1-\theta z)^p E(\sqrt{z}) = \frac{\pi}{2}  \sum_{n=0}^\infty u_n z^n, \quad |z| < \frac{1}{|\theta|}.
\end{equation*}
Then $u_0 = 1$, $u_1 = - (4 \theta  p+1) /4$, $u_2 = \frac{1}{64} \left(32 \theta ^2 p^2-32 \theta ^2 p+16 \theta  p-3\right)$ and
\begin{equation*}
u_{n+1} = \alpha_0(n) u_n + \alpha_1(n) u_{n-1} + \alpha_2(n) u_{n-2},
\end{equation*}
where
\begin{equation*}
\begin{aligned}
\alpha_0(n) &= \frac{(8 \theta +4) n^2-8 \theta  n p-4 \theta  p-1}{4 (n+1)^2},\\
\alpha_1(n) &= \frac{\theta  \left(-2 \theta  (-n+p+1)^2-(2 n-1) (2 n-2 p-3)\right)}{2 (n+1)^2},\\
\alpha_2(n) &= \frac{\theta ^2 (2 n-2 p-5) (2 n-2 p-3)}{4 (n+1)^2}.
\end{aligned}
\end{equation*}
\end{cor}

Taking $(a,b,c)=(1/2,1/2,1)$ and $(a,b,c)=(-1/2,1/2,1)$ in Theorem \ref{thm-e-arctan-F}, we obtain the following corollaries.

\begin{cor}
Let $p, z \in \mathbb{C}$. Suppose that
\begin{equation*}
e^{-p \arctan z} K(\sqrt{z}) = \frac{\pi}{2}  \sum_{n=0}^\infty u_n z^n, \quad |z|<1.
\end{equation*}
Then $u_0 = 1$, $u_1 = \frac{1}{4} (1-4 p)$, $u_2 = \frac{1}{64} \left(32 p^2-16 p+9\right)$, $u_3 = \frac{1}{768} \left(-128 p^3+96 p^2+148 p+75\right)$,
\begin{equation*}
u_4 = \frac{\left(2048 p^4-2048 p^3-12928 p^2-704 p+3675\right)}{49152}
\end{equation*}
and
\begin{equation*}
u_{n+1} = \beta_0(n) u_n + \beta_1(n) u_{n-1} + \beta_2(n) u_{n-2} + \beta_3(n) u_{n-3} + \beta_4(n) u_{n-4},
\end{equation*}
where
\begin{equation*}
\begin{aligned}
\beta_0(n) &= \frac{(2 n+1) (2 n-4 p+1)}{4 (n+1)^2},\\
\beta_1(n) &= -\frac{2 n^2-2 n (p+2)+p^2+2}{(n+1)^2},\\
\beta_2(n) &= \frac{4 n^2-4 n (p+3)+2 p (p+5)+9}{2 (n+1)^2}, \\
\beta_3(n) &= -\frac{(n-3) (n-2 p-3)}{(n+1)^2}, \\
\beta_4(n) &= \frac{(7-2 n)^2}{4 (n+1)^2}.
\end{aligned}
\end{equation*}
\end{cor}

\begin{cor}
Let $p, z \in \mathbb{C}$. Suppose that
\begin{equation*}
e^{-p \arctan z} E(\sqrt{z}) = \frac{\pi}{2}  \sum_{n=0}^\infty u_n z^n, \quad |z|<1.
\end{equation*}
Then $u_0 = 1$, $u_1 = \frac{1}{4} (-4 p-1)$, $u_2 = \frac{1}{64} \left(32 p^2+16 p-3\right)$, $u_3 = \frac{1}{768} \left(-128 p^3-96 p^2+292 p-15\right)$,
\begin{equation*}
u_4 = \frac{\left(2048 p^4+2048 p^3-17536 p^2-3136 p-525\right)}{49152}
\end{equation*}
and
\begin{equation*}
u_{n+1} = \beta_0(n) u_n + \beta_1(n) u_{n-1} + \beta_2(n) u_{n-2} + \beta_3(n) u_{n-3} + \beta_4(n) u_{n-4},
\end{equation*}
where
\begin{equation*}
\begin{aligned}
\beta_0(n) &= \frac{(2 n+1) (2 n-4 p-1)}{4 (n+1)^2},\\
\beta_1(n) &= -\frac{2 n^2-2 n (p+2)+p^2+p+2}{(n+1)^2},\\
\beta_2(n) &= \frac{4 n^2-4 n (p+4)+2 p (p+5)+15}{2 (n+1)^2}, \\
\beta_3(n) &= -\frac{n^2-2 n (p+3)+7 p+9}{(n+1)^2}, \\
\beta_4(n) &= \frac{(2 n-9) (2 n-7)}{4 (n+1)^2}.
\end{aligned}
\end{equation*}
\end{cor}

Taking $(a,b,c)=(1/2,1/2,1)$ and $(a,b,c)=(-1/2,1/2,1)$ in Theorems \ref{thm-sin-F} and \ref{thm-cos-F}, we obtain the following corollaries.

\begin{cor}
Let $p, z \in \mathbb{C}$. Suppose that
\begin{equation*}
\sin(pz) K(\sqrt{z}) = \sum_{n=0}^\infty u_n z^n, \quad \cos(pz) K(\sqrt{z}) = \sum_{n=0}^\infty v_n z^n, \quad |z|<1.
\end{equation*}
Then both $u_n$ and $v_n$ satisfies
\begin{equation*}
u_{n+1} = \sum_{i=0}^9 \gamma_i(n) u_{n-i}, \quad n \geq 9,
\end{equation*}
where
\begin{equation*}
\gamma_0(n) = \frac{3 (n-1) n+1}{(n+1)^2},
\end{equation*}
\begin{equation*}
\gamma_1(n) = \frac{4 n \left(2 n \left(8 n \left((n-8) p^2-n+5\right)+172 p^2-79\right)-392 p^2+145\right)+608 p^2-199}{16 n^2 (n+1)^2},
\end{equation*}
\begin{equation*}
\gamma_2(n) = \frac{(5-2 n)^2 (12 (n-4) n+49)-16 \left(n \left(4 n \left(4 n^2-46 n+191\right)-1381\right)+911\right) p^2}{16 n^2 (n+1)^2},
\end{equation*}
\begin{equation*}
\gamma_3(n) = \frac{-\left(4 n^2-24 n+35\right)^2+16 (8 (n-6) n+67) p^4+4 (8 n (3 n (4 (n-15) n+331)-2405)+17271) p^2}{16 n^2 (n+1)^2},
\end{equation*}
\begin{equation*}
\gamma_4(n) = \frac{p^2 \left(2 n \left(-8 n \left(n (2 n-37)+4 p^2+253\right)+248 p^2+6089\right)-934 p^2-13627\right)}{2 n^2 (n+1)^2},
\end{equation*}
\begin{equation*}
\gamma_5(n) = \frac{p^2 \left(96 n (2 n-19) p^2+8 n (2 n ((n-22) n+179)-1281)+16 p^4+4264 p^2+13627\right)}{4 n^2 (n+1)^2},
\end{equation*}
\begin{equation*}
\gamma_6(n) = -\frac{p^4 \left(8 n (4 n-45)+16 p^2+999\right)}{n^2 (n+1)^2},
\end{equation*}
\begin{equation*}
\gamma_7(n) = \frac{p^4 \left(8 (n-13) n+24 p^2+333\right)}{n^2 (n+1)^2},
\end{equation*}
\begin{equation*}
\gamma_8(n) = -\frac{16 p^6}{n^2 (n+1)^2},
\end{equation*}
\begin{equation*}
\gamma_9(n) = \frac{4 p^6}{n^2 (n+1)^2}.
\end{equation*}
\end{cor}

\begin{cor}
Let $p, z \in \mathbb{C}$. Suppose that
\begin{equation*}
\sin(pz) E(\sqrt{z}) = \sum_{n=0}^\infty u_n z^n, \quad \cos(pz) E(\sqrt{z}) = \sum_{n=0}^\infty v_n z^n, \quad |z|<1.
\end{equation*}
Then both $u_n$ and $v_n$ satisfies
\begin{equation*}
u_{n+1} = \sum_{i=0}^9 \gamma_i(n) u_{n-i}, \quad n \geq 9,
\end{equation*}
where
\begin{equation*}
\gamma_0(n) = \frac{n (3 n-5)+1}{(n+1)^2},
\end{equation*}
\begin{equation*}
\gamma_1(n) = \frac{4 n \left(4 n \left(n \left(4 (n-8) p^2-3 n+16\right)+86 p^2-29\right)-392 p^2+87\right)+608 p^2-99}{16 n^2 (n+1)^2},
\end{equation*}
\begin{equation*}
\gamma_2(n) = \frac{(3-2 n)^2 (2 n-7) (2 n-5)-16 (n (8 n (2 (n-12) n+103)-1523)+1021) p^2}{16 n^2 (n+1)^2},
\end{equation*}
\begin{equation*}
\gamma_3(n) = \frac{p^2 \left(16 n \left(6 n^3-96 n^2+2 (n-6) p^2+561 n-1426\right)+268 p^2+21305\right)}{4 n^2 (n+1)^2},
\end{equation*}
\begin{equation*}
\gamma_4(n) = \frac{p^2 \left(2 n \left(-8 n \left(2 (n-20) n+4 p^2+295\right)+256 p^2+7615\right)-3 \left(330 p^2+6049\right)\right)}{2 n^2 (n+1)^2},
\end{equation*}
\begin{equation*}
\gamma_5(n) = \frac{p^2 \left(8 \left(24 (n-10) n p^2+(n-12) n (2 (n-12) n+139)+2 p^4\right)+9 \left(524 p^2+2145\right)\right)}{4 n^2 (n+1)^2},
\end{equation*}
\begin{equation*}
\gamma_6(n) = -\frac{p^4 \left(32 (n-12) n+16 p^2+1141\right)}{n^2 (n+1)^2},
\end{equation*}
\begin{equation*}
\gamma_7(n) = \frac{2 p^4 \left(4 (n-14) n+3 \left(4 p^2+65\right)\right)}{n^2 (n+1)^2},
\end{equation*}
\begin{equation*}
\gamma_8(n) = -\frac{16 p^6}{n^2 (n+1)^2},
\end{equation*}
\begin{equation*}
\gamma_9(n) = \frac{4 p^6}{n^2 (n+1)^2}.
\end{equation*}
\end{cor}

Taking $(a,b,c)=(1/2,1/2,1)$ and $(a,b,c)=(-1/2,1/2,1)$ in Theorems \ref{thm-sinh-F} and \ref{thm-cosh-F}, we obtain the following corollaries.

\begin{cor}
Let $p, z \in \mathbb{C}$. Suppose that
\begin{equation*}
\sinh(pz) K(\sqrt{z}) = \sum_{n=0}^\infty u_n z^n, \quad \cosh(pz) K(\sqrt{z}) = \sum_{n=0}^\infty v_n z^n, \quad |z|<1.
\end{equation*}
Then both $u_n$ and $v_n$ satisfies
\begin{equation*}
u_{n+1} = \sum_{i=0}^9 \delta_i(n) u_{n-i}, \quad n \geq 9,
\end{equation*}
where
\begin{equation*}
\delta_0(n) = \frac{3 (n-1) n+1}{(n+1)^2},
\end{equation*}
\begin{equation*}
\delta_1(n) = \frac{4 n \left(2 n \left(-8 n \left((n-8) p^2+n-5\right)-172 p^2-79\right)+392 p^2+145\right)-608 p^2-199}{16 n^2 (n+1)^2},
\end{equation*}
\begin{equation*}
\delta_2(n) = \frac{16 \left(n \left(4 n \left(4 n^2-46 n+191\right)-1381\right)+911\right) p^2+(12 (n-4) n+49) (5-2 n)^2}{16 n^2 (n+1)^2},
\end{equation*}
\begin{equation*}
\delta_3(n) = -\frac{\left(4 n^2-24 n+35\right)^2-16 (8 (n-6) n+67) p^4+4 (8 n (3 n (4 (n-15) n+331)-2405)+17271) p^2}{16 n^2 (n+1)^2},
\end{equation*}
\begin{equation*}
\delta_4(n) = \frac{p^2 \left(2 n \left(8 n \left(n (2 n-37)-4 p^2+253\right)+248 p^2-6089\right)-934 p^2+13627\right)}{2 n^2 (n+1)^2},
\end{equation*}
\begin{equation*}
\delta_5(n) = \frac{p^2 \left(96 n (2 n-19) p^2-8 n (2 n ((n-22) n+179)-1281)-16 p^4+4264 p^2-13627\right)}{4 n^2 (n+1)^2},
\end{equation*}
\begin{equation*}
\delta_6(n) = \frac{p^4 \left(8 n (45-4 n)+16 p^2-999\right)}{n^2 (n+1)^2},
\end{equation*}
\begin{equation*}
\delta_7(n) = \frac{p^4 \left(8 (n-13) n-24 p^2+333\right)}{n^2 (n+1)^2},
\end{equation*}
\begin{equation*}
\delta_8(n) = \frac{16 p^6}{n^2 (n+1)^2},
\end{equation*}
\begin{equation*}
\delta_9(n) = -\frac{4 p^6}{n^2 (n+1)^2}.
\end{equation*}
\end{cor}

\begin{cor}
Let $p, z \in \mathbb{C}$. Suppose that
\begin{equation*}
\sinh(pz) E(\sqrt{z}) = \sum_{n=0}^\infty u_n z^n, \quad \cosh(pz) E(\sqrt{z}) = \sum_{n=0}^\infty v_n z^n, \quad |z|<1.
\end{equation*}
Then both $u_n$ and $v_n$ satisfies
\begin{equation*}
u_{n+1} = \sum_{i=0}^9 \delta_i(n) u_{n-i}, \quad n \geq 9,
\end{equation*}
where
\begin{equation*}
\delta_0(n) = \frac{n (3 n-5)+1}{(n+1)^2},
\end{equation*}
\begin{equation*}
\delta_1(n) = \frac{4 n \left(4 n \left(-3 n^2-2 (2 (n-8) n+43) p^2+16 n-29\right)+392 p^2+87\right)-608 p^2-99}{16 n^2 (n+1)^2},
\end{equation*}
\begin{equation*}
\delta_2(n) = \frac{16 (n (8 n (2 (n-12) n+103)-1523)+1021) p^2+(2 n-7) (2 n-5) (3-2 n)^2}{16 n^2 (n+1)^2},
\end{equation*}
\begin{equation*}
\delta_3(n) = \frac{p^2 \left(16 n \left(-6 n^3+96 n^2+2 (n-6) p^2-561 n+1426\right)+268 p^2-21305\right)}{4 n^2 (n+1)^2},
\end{equation*}
\begin{equation*}
\delta_4(n) = \frac{p^2 \left(2 n \left(8 n \left(2 (n-20) n-4 p^2+295\right)+256 p^2-7615\right)-990 p^2+18147\right)}{2 n^2 (n+1)^2},
\end{equation*}
\begin{equation*}
\delta_5(n) = \frac{p^2 \left(192 (n-10) n p^2-8 (n-12) n (2 (n-12) n+139)-16 p^4+9 \left(524 p^2-2145\right)\right)}{4 n^2 (n+1)^2},
\end{equation*}
\begin{equation*}
\delta_6(n) = \frac{p^4 \left(-32 (n-12) n+16 p^2-1141\right)}{n^2 (n+1)^2},
\end{equation*}
\begin{equation*}
\delta_7(n) = \frac{2 p^4 \left(4 (n-14) n-12 p^2+195\right)}{n^2 (n+1)^2},
\end{equation*}
\begin{equation*}
\delta_8(n) = \frac{16 p^6}{n^2 (n+1)^2},
\end{equation*}
\begin{equation*}
\delta_9(n) = -\frac{4 p^6}{n^2 (n+1)^2}.
\end{equation*}
\end{cor}

\end{document}